\documentclass[preprint,aos]{imsart}

\usepackage{amsmath,amssymb,amsthm,amsfonts,amstext,amsbsy,amscd}

\usepackage{comment}
\usepackage{graphicx}
\usepackage{mathrsfs}
\usepackage{color}
\usepackage[authoryear]{natbib}
\RequirePackage[colorlinks,citecolor=blue,urlcolor=blue]{hyperref}

\numberwithin{equation}{section}

\swapnumbers

\newtheorem{satz}{Satz}[section]

\newtheorem{theorem}[satz]{Theorem}
\newtheorem{proposition}[satz]{Proposition}
\newtheorem{corollary}[satz]{Corollary}
\newtheorem{lemma}[satz]{Lemma}
\newtheorem{assumption}[satz]{Assumption}

\newtheorem{definition}[satz]{Definition}
\theoremstyle{remark}

\newtheorem{remark}[satz]{Remark}

\newtheorem{example}[satz]{Example}

\DeclareMathOperator{\E}{{\mathbb E}}

\DeclareMathOperator{\R}{{\mathbb R}}
\DeclareMathOperator{\C}{{\mathbb C}}

\DeclareMathOperator{\N}{{\mathbb N}}

\DeclareMathOperator{\HH}{{\mathbb H}}

\DeclareMathOperator{\PP}{{\mathbb P}}

\DeclareMathOperator{\spann}{span}
\DeclareMathOperator{\trace}{trace}

  \DeclareMathOperator{\rank}{rank}
\DeclareMathOperator{\esssup}{esssup}
\DeclareMathOperator{\diag}{diag} 
 \DeclareMathOperator{\Id}{Id}
\DeclareMathOperator{\argmin}{argmin}
\DeclareMathOperator{\Var}{Var}

\DeclareMathOperator{\vek}{vec}

\providecommand{\eps}{\varepsilon}
\renewcommand{\phi}{\varphi}
\renewcommand{\theta}{\vartheta}
\renewcommand{\subset}{\subseteq}

\renewcommand{\cdot}{{\scriptstyle \bullet} }
\providecommand{\abs}[1]{\lvert #1 \rvert}
\providecommand{\norm}[1]{\lVert #1 \rVert}
\providecommand{\bnorm}[1]{{\Bigl\lVert #1 \Bigr\rVert}}

\providecommand{\scapro}[2]{\langle #1,#2 \rangle}
\providecommand{\floor}[1]{\lfloor #1 \rfloor}

\providecommand{\ceil}[1]{\lceil #1 \rceil}

\renewcommand{\le}{\leqslant}
\renewcommand{\ge}{\geqslant}
\providecommand{\COV}{\C\!\operatorname{ov}}

\providecommand{\mr}{\color{blue}}

\allowdisplaybreaks[4]

\begin{document}

\begin{frontmatter}
\title{Inference on the maximal rank of\\ time-varying covariance matrices\\ using high-frequency data}
\runtitle{Rank inference for time-varying covariance matrices}

\begin{aug}
\author[A]{\fnms{Markus} \snm{Rei\ss}\ead[label=e1]{mreiss@math.hu-berlin.de}}
\and
\author[B]{\fnms{Lars} \snm{Winkelmann}\ead[label=e2]{lars.winkelmann@fu-berlin.de}}

\address[A]{Institute of Mathematics, Humboldt-Universit\"at zu Berlin,
\printead{e1}}

\address[B]{School of Business and Economics, Freie Universit\"at Berlin,
\printead{e2}}
\end{aug}

\begin{abstract}
We study the rank of the instantaneous or spot covariance matrix $\Sigma_X(t)$ of a multidimensional continuous semi-martingale $X(t)$. Given high-frequency observations $X(i/n)$, $i=0,\ldots,n$, we test the null hypothesis $\rank(\Sigma_X(t))\le r$ for all $t$ against local alternatives where the average $(r+1)$st eigenvalue is larger than some signal detection rate $v_n$. 

A major problem is that the inherent averaging in local covariance statistics produces a bias that distorts the rank statistics. We show that the bias depends on the regularity and a spectral gap of $\Sigma_X(t)$. We establish explicit matrix perturbation and concentration results that provide non-asymptotic uniform critical values and optimal signal detection rates $v_n$. This leads to a rank estimation method via sequential testing.  For a class of stochastic volatility models, we determine data-driven critical values via normed p-variations of estimated local covariance matrices. The methods are illustrated by simulations and an application to high-frequency data of U.S. government bonds.

\end{abstract}



\end{frontmatter}

\section{Introduction}
Consider a continuous-time $\R^d$-valued process $(X(t),t\in[0,1])$ with  instantaneous or spot covariance matrix
$\Sigma_X(t)=\lim_{h\to  0}h^{-1}\COV(X(t+h)-X(t))\in\R^{d\times d}.$
We address the problem of testing and estimating the maximal rank of $\Sigma_X(t)$ on $t\in[0,1]$, using discrete observations $X(i/n)$, $i=0,\ldots,n$. These high-frequency observations allow for inference without specific modeling assumptions on $\Sigma_X$.

In econometrics, the maximal rank $r$ may be the number of risk factors involved in financial asset prices $X$ or the dimension of the state space of a continuous-time factor model. See \citet{jacod2013} for a discussion of market completeness and \citet{ait2017} who provide a link between continuous-time factor models and principal component analysis (PCA). $\Sigma_X(t)$ is well defined if $X(t)$ is modeled by an $L^2$-semi-martingale.  As discussed below,
the basic case is given by a Brownian martingale
\begin{equation}\label{EqX}
 dX(t)=\sigma_X(t)dB(t),\quad t\in[0,1],
 \end{equation}
with an $m$-dimensional standard Brownian motion $B(t)$, $\sigma_X(t)\in\R^{d\times m}$ and  $\Sigma_X(t)=\sigma_X(t)\sigma_X(t)^\top\in\R^{d\times d}$. For deterministic $\Sigma_X$, we aim at inferring the rank of $\Sigma_X(t)$, $t\in [0,1]$, from the observed independent increments
\begin{equation}\label{EqXincr}
X(\tfrac in)-X(\tfrac{i-1}n)\sim N\Big(0,\int_{(i-1)/n}^{i/n}\Sigma_X(t)\,dt\Big),\quad i=1,\ldots,n.
\end{equation}
This formulation poses a fundamental problem in nonparametric statistics, whose analysis will turn out to be non-standard.

A natural statistic for the integrated covariance on a block $I_k=[kh,(k+1)h]\subset[0,1]$ for some small $h>0$ with $nh\in\N$ is the local  empirical ({\it realized}) covariance matrix
\begin{equation}\label{EqSigmaXhat}
\hat\Sigma_X^{kh}:=\frac1h\sum_{i=1}^{nh}\big(X(kh+\tfrac in)-X(kh+\tfrac{i-1}n)\big)\big(X(kh+\tfrac in)-X(kh+\tfrac{i-1}n)\big)^\top.
\end{equation}
The realized covariance matrix estimates the average covariance matrix
\begin{equation}\Sigma_X^{kh}:=\E\big[\hat\Sigma_X^{kh}\big]=\frac1h\int_{I_k}\Sigma_X(t)\,dt\label{EqESigmakh}
\end{equation}
with an error of size ${\mathcal O}_P((nh)^{-1/2})$.
The average $\Sigma_X^{kh}$ can have full rank even though $\rank(\Sigma_X(t))\le r < d$ holds for all $t\in I_k$. The potential inflation in the rank is a consequence of averaging spot covariance matrices with time-varying eigenspaces. We highlight the theoretical and empirical relevance of such time-varying features in deterministic and stochastic volatility models as well as financial data. This confirms \cite{jacod2008} who state that $\Sigma_X^{kh}$ gives no insight on the rank of $\Sigma_X(t)$. We show, however, while the $(r+1)$st eigenvalue $\lambda_{r+1}(\Sigma_X^{kh})$ may be non-zero, it is
still small for small block sizes $h$ when $\Sigma_X(t)$ varies smoothly. Precise perturbation bounds for matrix averages indicate that the maximal size of $\lambda_{r+1}(\Sigma_X^{kh})$ does not only depend on the regularity $\beta$ of $\Sigma_X$, but also significantly on the existence of a spectral gap $\underline\lambda_r$. Small  eigenvalues $\lambda_{r+1}(\hat\Sigma_X^{kh})$ across blocks $I_k$ should favour the acceptance of the null hypothesis, while large values should lead to a rejection. This paper develops such a test with specific attention to the block size $h$ and to non-asymptotic critical values.  At an abstract level a classical bias-variance dilemma seems to dictate the choice of the block size $h$ in \eqref{EqSigmaXhat}.

By studying the null hypothesis $\lambda_{r+1}(\Sigma_X(t))=0$ for all $t$ versus the alternative $\lambda_{r+1}(\Sigma_X(t))>0$ for some $t$, we embed our rank test into a signal detection framework.
 In particular, we allow for local alternatives where  $\lambda_{r+1}(\Sigma_X(t))$ has average size at least $v_n$ where the {\it signal strength} $v_n$ tends to zero as $n\to\infty$. In contrast to simple consistency results under a fixed alternative, this approach not only reveals the approximate deviations from $\lambda_{r+1}(\Sigma_X(t))=0$ if the test accepts, but also allows to establish the minimax optimal signal detection rate $v_n=\min(n^{-\beta},\underline\lambda_r^{-1}n^{-2\beta})$ in the sense of \citet{ingster2012}. This rate is
 attained by considering the average of the empirical eigenvalues $\lambda_{r+1}(\hat\Sigma_X^{kh})$ over all blocks and choosing a block size $h$ of order $n^{-1}$. The bias-variance dilemma disappears at the level of rates because of the heteroskedastic error of $\hat\Sigma^{kh}$, which is natural for variance-type estimators, and the small size of $\lambda_{r+1}(\Sigma_X^{kh})$ under the null.
 The heteroskedasticity in combination with the bias bound also explains the surprisingly fast detection rate $v_n$ compared to the classical estimation rates $n^{-\beta/(2\beta+1)}$ and $n^{-1/2}$ for the spot covariance $\Sigma_X(t)$ and the integrated covariance $\int_0^1\Sigma_X(t)dt$, respectively.

 For standard stochastic volatility models with regularity $\beta=1/2$ the rank detection rate is $v_n=\min(n^{-1/2},(\underline\lambda_rn)^{-1})$, opposed to the rates $n^{-1/4}$ and $n^{-1/2}$  for the spot covariance and integrated covariance estimation, see e.g. \citet{FanWang2008}.
Compared to classical signal detection the roles of hypothesis and alternative are in a certain sense reversed. We shall understand this by an underlying concave functional instead of a convex constraint for classical regression, see e.g. \citet{juditsky2002} for  the interplay between regularity and convexity constraints.

The critical values of our tests depend on a regularity bound for the covariance matrix function $\Sigma_X$ under the null hypothesis. This may be available from previous observations or other side information, but with models for $\Sigma_X$ at hand also data-driven critical values are feasible. We present an approach assuming a semi-martingale model for $\Sigma_X$ itself. The regularity in a Besov scale is $\beta=1/2$, while the corresponding constant (similar to a Besov or H\"older norm) can be estimated by a $p$-variation of norms of empirical covariance matrix increments. This setting is similar to the scalar case in  \cite{Vetter2015}, who estimates the volatility of volatility, but it involves non-differentiable matrix norms. A feasible central limit theorem 
permits a data-driven calibration of  critical values under stochastic volatility and yields asymptotically the correct size. 

On the technical side we need to control the sum of empirical eigenvalues $\lambda_{r+1}(\hat\Sigma_X^{kh})$ over blocks $I_k$. $\hat\Sigma_X^{kh}$ in \eqref{EqSigmaXhat} is a convex combination of independent Wishart matrices with different population covariance matrices. To obtain non-asymptotic and uniform critical values, we refine general matrix Bernstein inequalities by \citet{tropp2012} and use Gaussian concentration. For the power analysis the standard lower-tail inequalities do not apply. Instead, we determine the specific eigenvalue density by establishing a stochastic domination property with respect to a fixed population covariance and use an entropy argument. Simple approaches are not feasible because the different population covariance matrices do not commute and are not sufficiently close to each other. Beyond the purposes of this paper, these results might find some independent interest.

Our method can be extended to cover idiosyncratic components, where we observe a process $Y$ satisfying
\begin{equation}\label{EqY}
 dY(t)=dX(t)+dZ(t),\text{ where } dZ(t)=\sigma_Z(t)\,dB_Z(t),
 \end{equation}
 with an independent $d$-dimensional Brownian motion $B_Z$, $\sigma_Z(t)\in\R^{d\times d}$ and spot covariances $\Sigma_Z(t)=\sigma_Z(t)\sigma_Z(t)^\top$ such that $\Sigma_Y(t)=\Sigma_X(t)+\Sigma_Z(t)$. Our methods directly translate to cases where a small upper bound for $\norm{\Sigma_Z(t)}$ is known, which gives a desirable robustness property.
 This provides a link to the wide class of factor models. In contrast to classical factor models for volatility \citep{LiYao2016,ait2017}, we allow the eigenspaces of $\Sigma_X(t)$ ({\it factor loading spaces}) to vary over time. 

Estimating and testing the rank of a covariance matrix has attracted a lot of attention under different angles. The case where the covariance $\Sigma_X$ in \eqref{EqXincr} is constant (or more generally has constant null space in time) is usually trivial because the rank of the empirical covariance matrix  equals the rank of $\Sigma_X$ as soon as $n\ge d$. Observing $Y(i/n)$, $i=0,\ldots,n$, with constant $\Sigma_X$ and $\Sigma_Z$ in \eqref{EqY} leads to a spiked covariance model \citep{johnstone2001}. In this i.i.d. framework and for high dimensions \cite{onatski2014} develop asymptotic tests on the presence of spikes  and \citet{cai2015} establish optimal rates for rank detection. In these works the spectral gap plays the role of a signal strength, while in our setting it controls the perturbations due to averaging  and leads to different exponents in the detection rate. We focus on the effect of time-varying covariances and do not study additional dimension asymptotics or nuclear norm penalisations as e.g. in \cite{christensen2021}. Nevertheless, our non-asymptotic results show the dependence on the dimension $d$ and the rank $r$  explicitly. Time-varying covariances appear naturally in time series analysis.  Inference on their rank, however, requires much more specific modeling assumptions on the discrete-time dynamics of $\Sigma_X$. We refer to \citet{LamYao2012} for an analysis in a stationary high-dimensional framework, to \citet{Su2017} for a local PCA approach to testing constant factor loadings over time and generally to the references therein for many further aspects.

The setting in our paper is close to that of \citet{jacod2013}, who consider  rank estimation and testing problems in a joint semi-martingale setup for $X(t)$ and $\sigma_X(t)$. They establish stable central limit theorems with convergence rate $n^{-1/2}$ for statistics that involve an additional randomisation step. Since the integer-valued rank is their target, there is no clear concept of convergence rates or local alternatives as in our signal detection framework. While they treat general joint semi-martingale models, we focus on the pathwise properties of $\Sigma_X$ in terms of regularity and spectral gap.   \citet{ait2017} detect the rank of the integrated covariance matrix $\int_0^1\Sigma_X(t)dt$ in a sparse high-dimensional setting and apply this to determine the number of factors given constant factor loading spaces. More recently, \citet{ait2019}  argue that the rank of the spot covariance matrix is often much more informative and provide empirical evidence based on the S\&P 100 index. They develop an asymptotic theory for the so-called realised eigenvalues $\int_0^t\lambda_j(\Sigma(s))\,ds$ and more general spectral functions, but they focus on the case of  covariance matrices with full rank. An inspiring  characterisation of the realised eigenvalue $\int_0^1\lambda_{r+1}(\Sigma(s))\,ds$ is given by \citet{jacod2008} as a natural distance to all semi-martingales with spot covariance matrix of rank at most $r$.

In Section \ref{SecSetting} we introduce the test statistics and discuss in theory and examples the impact of averaging on the eigenvalue perturbation. Our data example refers to U.S.\ government bonds and the term structure of interest rates.
The main results  are presented in Section \ref{SecResults}. This includes the analysis of the tests under null hypothesis and local alternatives, the derivation of the optimal signal detection rate as well as properties of a corresponding rank estimation method. Section \ref{SecStochVol} studies the specific case of stochastic volatility models and  provides data-driven critical values based on a matrix norm $p$-variation. These techniques are applied to inference on the rank in the bond data. Mathematical tools and more technical proofs are delegated to the appendix. Some of the matrix deviation and concentration results there are stated in wider generality
because they might prove useful in other circumstances as well.

\section{\bf Setting, examples and first results} \label{SecSetting}

Let us first fix some notation. We write $a_n\lesssim b_n$ or $a_n={\mathcal O}(b_n)$ if $a_n\le Cb_n$ for some constant $C>0$ and all $n$. By $a_n\thicksim b_n$ we mean  $a_n\lesssim b_n$ and $b_n\lesssim a_n$.  With ${\mathcal O}_P(b_n)$ and ${o}_P(b_n)$ we denote random variables $X_n$ (also matrix-valued) such that $(b_n^{-1}X_n)_{n\ge 1}$ remain bounded in probability and tend to zero in probability, respectively. Similarly, ${\mathcal O}_{L^p}(b_n)$ stands for random variables $X_n$ with $\E[\abs{X_n}^p]^{1/p}\lesssim b_n$.  For two random variables $X,Y$ with the same law we write $X\stackrel{d}{=}Y$.

The canonical basis in $\R^d$ is denoted by $e_1,\ldots,e_d$, $E_{i,j}=e_ie_j^\top\in\R^{d\times d}$ is an elementary matrix and $I_d$ denotes the identity matrix in $\R^{d\times d}$.
We introduce the sets of  matrices
\[\R^{d\times d}_{sym}:=\{S\in\R^{d\times d}\,|\, S=S^\top\},\quad \R^{d\times d}_{spd}:=\{S\in\R_{sym}^{d\times d}\,|\, S\text{ positive semi-definite}\}.\] For $A,B\in\R_{sym}^{d\times d}$ the partial order $A\le B$ says that $B-A\in \R_{spd}^{d\times d}$.  For $S\in\R_{sym}^{d\times d}$ we consider the ordered eigenvalues $\lambda_{max}(S)=\lambda_1(S)\ge\cdots\ge\lambda_d(S)=\lambda_{min}(S)$ (according to their multiplicities). For matrices $S,T\in\R^{d\times d}$ the Hilbert-Schmidt (or Frobenius) scalar product is $\scapro{S}{T}_{HS}=\trace(T^\top S)^{1/2}$   and the spectral norm is given by $\norm{S}=\max_{v\in\R^d,\norm{v}=1}\norm{Sv}$.

We split the interval $[0,1]$ into blocks $I_k=[kh,(k+1)h]$ for some block length $h$ with $nh,h^{-1}\in\N$ and $k=0,\ldots,h^{-1}-1$. On each block $I_k$ we consider the empirical or realised covariance $\hat\Sigma_X^{kh}$ from \eqref{EqSigmaXhat} and its mean $\Sigma_X^{kh}$ in \eqref{EqESigmakh}. Assuming  model \eqref{EqX} for $X$ with deterministic or independent $\sigma_X(t)$, we obtain from \eqref{EqXincr} (conditionally on $\Sigma_X$)
\begin{align}
\COV\big(\vek(\hat\Sigma_X^{kh})\big)&=\frac1{nh^2}\sum_{i=1}^{nh}\Big(\int_{kh+(i-1)/n}^{kh+i/n}\Sigma_X(t)\,dt\Big)^{\otimes 2}{\mathcal Z}_d,\label{EqCovSigmakh}
\end{align}
where we follow \citet{magnus1979} and vectorize matrices $A$ by $\vek(A)$, employ the Kronecker product $\otimes$ for matrices and use the matrix ${\mathcal Z_d}=\COV(\vek(ZZ^\top))\in\R^{d^2\times d^2}$ for $Z\sim N(0,I_d)$.
If the process $Y$ in \eqref{EqY} is observed, the corresponding covariance matrices in terms of $Y$ are denoted by $\hat\Sigma_Y^{kh}$ and $\Sigma_Y^{kh}$. In this section we focus on $X$.

We  aim at a level-$\alpha$ test $\phi_\alpha$ for the null hypothesis $\max_{t\in[0,1]}\rank(\Sigma_X(t))\le r$ of the form
\begin{equation}\label{Eqphi0} \phi_\alpha:={\bf 1}\big(T_{n,h}>\kappa_\alpha\big)\text{ with } T_{n,h}:=\sum_{k=0}^{h^{-1}-1} h\lambda_{r+1}\big(\hat\Sigma_X^{kh}\big).
\end{equation}
Observe that the number $nh$ of observations on each block should be at least $r+1$ because otherwise $\lambda_{r+1}\big(\hat\Sigma_X^{kh}\big)=0$ holds.
The choice of the block size $h$ and of the critical value $\kappa_\alpha$ requires a deeper understanding of eigenvalue deviations under averaging and stochastic errors.

The following is an easy consequence of the variational characterisation of eigenvalues, see e.g. \citet[Proposition 1.3.4]{tao2012}.
\begin{lemma}\label{Lemtracer}
The remaining partial trace map $S\mapsto \trace_{>r}(S)=\sum_{j=r+1}^d\lambda_j(S)$ is concave for $S\in\R^{d\times d}_{sym}$ and $0\le r<d$. In particular,
the remaining partial trace of an average is larger than the average over the remaining partial trace:
\begin{equation}\label{Eqtracer}
 \trace_{> r}\big(\Sigma_X^{kh}\big)=\trace_{> r}\Big(\frac1h\int_{I_k}\Sigma_X(t)dt\Big)\ge \frac1h\int_{I_k}\trace_{> r}(\Sigma_X(t))dt.
\end{equation}
Consequently, $\rank(\Sigma_X^{kh})\ge \esssup_{t\in I_k}\rank(\Sigma_X(t))$ holds.
  \end{lemma}

While Lemma \ref{Lemtracer} shows that averaging can only increase the rank, we want to understand how large  $\lambda_{r+1}(\Sigma_X^{kh})$ can become if $\rank(\Sigma_X(t))=r$ holds for all $t$. Three examples are discussed before we state a precise perturbation bound. The first bivariate example shows the eigenvalue perturbation by rotating eigenvectors under H\"older regularity. For $\beta\in(0,1]$ we use the matrix-valued $\beta$-H\"older ball of radius $L>0$
\[ C^\beta(L):=\{\Sigma: [0,1]\to\R_{spd}^{d\times d}\,|\, \norm{\Sigma(t)-\Sigma(s)}\le L\abs{t-s}^\beta\text{ for all } t,s\in [0,1]\}.
\]
The impact of a spectral gap is clearly exposed.

\begin{example}\label{ExRank}
 Let $d=2$, $r=1$, $\underline\lambda_1\ge  h^\beta/\sqrt 2$ and consider
\[ \Sigma_X(t)=\begin{pmatrix} \underline\lambda_1 & h^\beta \sin(2\pi t/h)\\ h^\beta \sin(2\pi t/h) & \underline\lambda_1^{-1}h^{2\beta}\sin^2(2\pi t/h)\end{pmatrix},\quad  \frac1h\int_0^h\Sigma_X(t)dt=\begin{pmatrix} \underline\lambda_1 & 0\\ 0 & (2\underline\lambda_1)^{-1}h^{2\beta}\end{pmatrix}.
\]
Then $\lambda_1(\Sigma_X(t))\ge \underline\lambda_1$ and $\lambda_2(\Sigma_X(t))=0$ hold for all $t\in[0,1]$.
We have $\Sigma_X\in C^\beta(L)$ for $L=6\pi$ (use
$\norm{\Sigma_X'(t)}
\le 6\pi h^{\beta-1}$). The average $h^{-1}\int_0^h\Sigma_X(t)dt$, however,
has rank 2 with
\[\lambda_{2}\Big(\frac1h\int_0^h\Sigma_X(t)dt\Big)=\tfrac12 \underline\lambda_1^{-1}h^{2\beta}.
\]
Note $\frac12 \underline\lambda_1^{-1}h^{2\beta}\le\underline\lambda_1$ by the assumption on $\underline\lambda_1$.
For $\underline\lambda_1\thicksim h^\beta$ the spectral gap $\lambda_1(\Sigma_X(t))-\lambda_2(\Sigma_X(t))$ has order $h^\beta$ and we obtain $\lambda_{2}(h^{-1}\int_0^h\Sigma_X(t)dt)\thicksim h^\beta$, which seems a natural deviation from $\lambda_2(\Sigma_X(t))=0$ for $\beta$-H\"older continuous $\Sigma$. Yet, a spectral gap of order $1$, i.e. $\underline\lambda_1\thicksim  1$, yields a much smaller quadratic deviation $\lambda_{2}(h^{-1}\int_0^h\Sigma_X(t)dt)\thicksim h^{2\beta}$ from zero.

In the extreme case $h=1$, PCA on the integrated covariance matrix $\int_0^1\Sigma_X(t)dt$ for $\underline\lambda_1=1$ would result in $33\%$ explained (population) variance by the second component although the spot covariances are all of rank one.
\end{example}

The  size of the eigenvalue perturbation in Example \ref{ExRank} is also attained in typical stochastic volatility models, which shows that this is not only  a worst case scenario.

\begin{example}\label{ExStochVol}
Stochastic volatility models are often based on $d$-dimensional Wishart processes  \citep{Bru1991}, compare the matrix square Ornstein-Uhlenbeck and affine processes. The Wishart process of dimension $d$ and index $r$ is given by $\Sigma_X(t)=\tilde B(t)^\top\tilde B(t)$, with an $(r\times d)$-dimensional  Brownian matrix $\tilde B(t)$ (i.e., the entries $\tilde B_{ij}(t)$ form independent Brownian motions). For general deterministic initial value $\tilde B(0)=b_0\in\R^{r\times d}$, $s_0=b_0^\top b_0\in\R^{d\times d}$ is the initial state of  $\Sigma_X(t)$.
The Wishart process forms a matrix-valued It\^o semi-martingale, satisfying the stochastic differential equation
\[ d\Sigma_X(t)=rI_ddt+\Sigma_X(t)^{1/2}d\check B(t)+\big(\Sigma_X(t)^{1/2}d\check B(t)\big)^\top,\quad \Sigma(0)=s_0,\]
with a $(d\times d)$-Brownian matrix $\check B(t)$, $t\ge 0$. A matrix square Ornstein-Uhlenbeck process, for example, satisfies the same equation, but with a more general drift term $(2\gamma\Sigma_X(t)+rI_d)$ for some back-driving parameter $\gamma<0$. On short time intervals, the drift term will become negligible so that  the Wishart process serves as a fundamental example.

For $r<d$ and an initialisation with $\rank(s_0)=r$ the Wishart process $\Sigma_X(t)$ has by definition rank at most $r$. Yet, the average $\frac1h\int_0^h\Sigma(t)dt$ has almost surely full rank with an $(r+1)$st eigenvalue of size ${\mathcal O}_P(h)$ for small $h>0$.

It is straight-forward to derive these properties in the case $d=2$, $r=1$, $\tilde B(t)=(1,0)+\breve B(t)$, with $\breve B(t)$ a planar Brownian motion starting in the origin. Then
\begin{align*}
\Sigma_X(t)&=\begin{pmatrix} (1+\breve B_1(t))^2 & (1+\breve B_1(t))\tilde B_2(t)\\ (1+\breve B_1(t))\breve B_2(t) & \breve B_2(t)^2\end{pmatrix},\\
\frac1h\int_0^h\Sigma_X(t)\,dt&=\begin{pmatrix} 1+{\mathcal O}_P(h^{1/2}) & \frac1h\int_0^h\breve B_2(t)dt+{\mathcal O}_P(h)\\ \frac1h\int_0^h\breve B_2(t)dt+{\mathcal O}_P(h) & \frac1h\int_0^h\breve B_2(t)^2dt\end{pmatrix}
\end{align*}
holds, implying $\lambda_2(\frac1h\int_0^h\Sigma_X(t)dt)=\int_0^h\breve B_2(t)^2dt-\frac1h(\int_0^h\breve B_2(t)dt)^2+{\mathcal O}_P(h^{3/2})$. By Brownian scaling we see the equality in law
\[\lambda_2\Big(\frac1h\int_0^h\Sigma_X(t)dt\Big)\stackrel{d}{=}h\Big(\int_0^1\breve B_2(t)^2dt-\Big(\int_0^1\breve B_2(t)dt\Big)^2\Big)+{\mathcal O}_P(h^{3/2})\]
and that the second eigenvalue is of size $h$. The spectral gap in $\Sigma_X(t)$ is $1+{\mathcal O}_P(h^{1/2})$ and the H\"older regularity $\beta$  is almost $1/2$, in line with the  order ${\mathcal O}(\underline\lambda_1^{-1}h^{2\beta})$ from Example \ref{ExRank}.

\begin{figure}[t]
\centering
\includegraphics[width=0.48\textwidth]{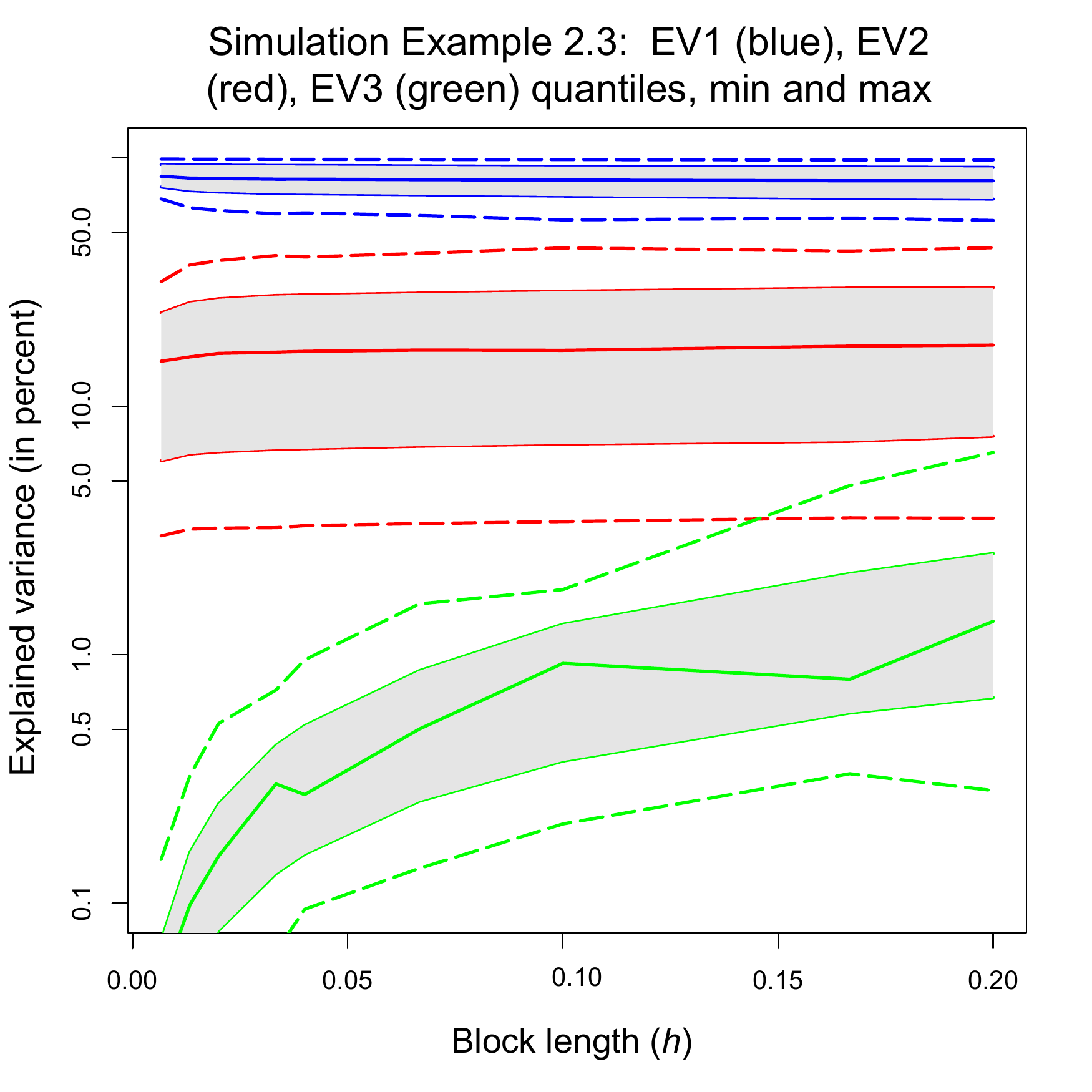} \includegraphics[width=0.48\textwidth]{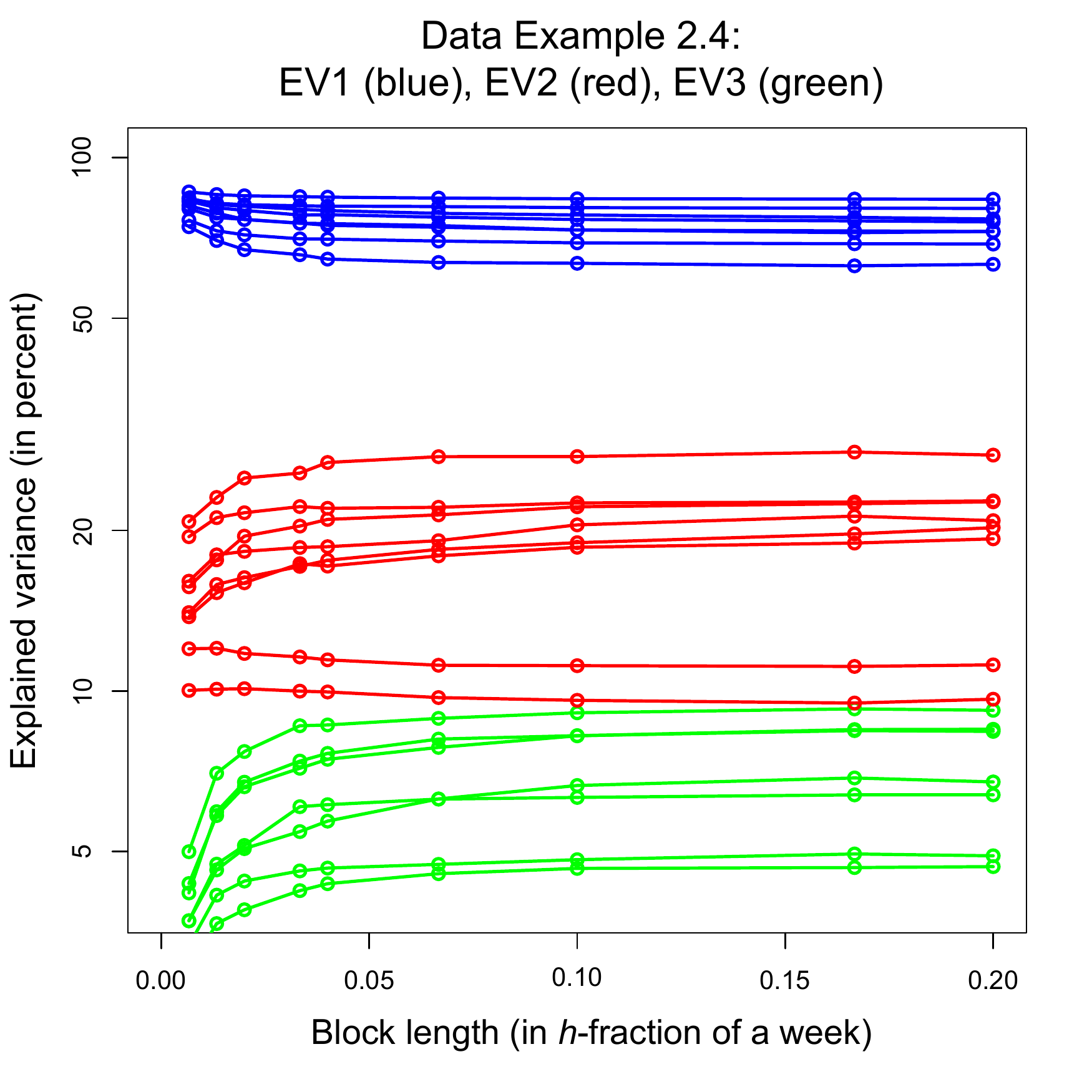}
 \caption{Dependence of eigenvalues of covariance matrices on the block length. Left:  simulated covariance matrices and their eigenvalues (EV) from Example \ref{ExStochVol}. Right:  eigenvalues of realised covariance matrices of U.S.\ government bonds in seven different week in 2020 from Example \ref{ExEmpirical}.   }\label{FigEx}
\end{figure}

For $d=3$, $r=2$, $n=1950$ and varying $h$ Figure \ref{FigEx}(left) shows statistics of the relative explained variance of component $j$ \begin{equation}\label{EqExplVar}
\frac{T_{n,h}^{(j)}}{T_{n,h}^{(1)}+T_{n,h}^{(2)}+T_{n,h}^{(3)}}\text{ with } T_{n,h}^{(j)}:=\sum_{k=0}^{h^{-1}-1}h\lambda_j(\hat\Sigma_X^{kh}),\quad j=1,2,3,
\end{equation}
on a logarithmic scale in 1000 Monte-Carlo simulations. Starting values of the two non-zero eigenvalues are 1 and 0.5. 
The range between the 10\% and 90\% quantiles for each $h$ is indicated by the gray-shaded areas. The simulation round that provides the median, min and max percentage for $h=0.2$ are given by the central solid and the two dashed lines. In relative terms the third eigenvalue (green) increases up to a maximal value of 5\% and a median of around 1\% if covariance matrices are estimated over blocks of length $h=0.2$, even though the third spot eigenvalues are all zero, which is attained in the limit $h\to 0$.
\end{example}

\begin{example}\label{ExEmpirical}
The popular \cite{NS87} model of the term structure of interest rates assumes a rank 2 covariance matrix across bonds of different maturities. The related two components are commonly known as level and slope factors and play an important role in asset pricing, risk management and macro-finance applications. Often larger models are proposed like the extension to three factors in \cite{DL06}.
We analyse $d=3$ U.S.\ government bond Exchange Traded Funds (ETFs): the iShares treasury of 3-7 years  (IEI), 7-10 years (IEF) and 10-20 years (TLH). The data is obtained from the Nasdaq through the data provider LOBSTER and includes one-minute data (6.5 trading hours per day) for the first 27 weeks in 2020. Using jump-truncated one-minute log returns over one week intervals ($n=1950$), Figure \ref{FigEx}(right) shows for seven selected weeks how the averaging of covariance matrices affects the  explained variances \eqref{EqExplVar} of the components. The third eigenvalue ($j=3$, green) that explains usually much less than 5\% of the total variance on blocks with 13 minutes length ($h=6.6\times10^{-3}$) becomes seemingly more important for larger blocks. Compared with daily data analysis (390 minutes, $h=0.2$), the magnitude roughly doubles.
We come back to the bond example in Section \ref{SecStochVol}, where Figure \ref{FigBondData} illustrates the significant dynamics of the empirical eigenvectors.
\end{example}
As a first step towards deriving the critical values $\kappa_\alpha$ in \eqref{Eqphi0}, we provide a general bound for eigenvalues of averages over low-rank matrices.  The proof in Appendix \ref{secMatrixConc} is inspired by the techniques in \citet{reiss2020}.

\begin{proposition}\label{PropPerturb}
Let $I\subset\R$ be an interval of length $\abs{I}$. Consider $S(t)\in \R^{d\times d}_{spd}$ with $\rank(S(t))\le r<d$ and $\lambda_r(S(t))\ge\underline\lambda_r\ge 0$ for $t\in I$. Then
\[ \lambda_{r+1}\Big(\frac1{\abs{I}}\int_I S(t)dt\Big) \le \big(\tfrac{2}{\underline\lambda_r}\Delta_2(S,I)^2\big)\wedge \Delta_1(S,I)\]
holds with the $L^p$-variations of $S$ on $I$
\[ \Delta_p(S,I):=\Big(\frac1{\abs{I}^2}\int_{I\times I} \norm{S(s)-S(t)}^p\,dt\,ds\Big)^{1/p},\quad 1\le p<\infty.
\]
\end{proposition}

\begin{remark}
As the proof reveals, asymptotically for $\Delta_2(S,I)\to 0$ and fixed $\underline\lambda_r>0$, the quadratic bound becomes $\frac1{2\underline\lambda_r}\Delta_2(S,I)^2$ since $\frac{\underline\lambda_r}{(\underline\lambda_r-\bar \lambda_{r+1})^2}\to 1$. This provides smaller asymptotic critical values.
\end{remark}

Notice that the bound of Proposition \ref{PropPerturb} is quadratic in case of a positive spectral gap $\underline\lambda_r$, which will later allow the detection of weaker signals (smaller eigenvalues).
For $\beta$-H\"older-continuous $S\in C^\beta(L)$, we have $\Delta_p(S,I)\le L\abs{I}^\beta$ and we obtain the order in Example \ref{ExRank}:
\[ \lambda_{r+1}\Big(\frac1{h}\int_0^h S(t)dt\Big) \lesssim \big(\tfrac{L^2}{\underline\lambda_r}h^{2\beta}\big)\wedge \big(Lh^\beta\big).\]
The bound in terms of an integrated variation criterion allows to weaken the H\"older regularity to an $L^p$-Besov regularity. This is important to cover stochastic volatility models as in Example \ref{ExStochVol} or in Section \ref{SecStochVol} below.

\section{\bf General results}\label{SecResults}

\subsection{\bf Behaviour under the null hypothesis}

Let us recall the definitions of Besov spaces $B^\beta_{p,\infty}$, see \cite{Cohen2003} for a nice survey. For $S\in L^p([0,1];\R^{d\times d})$ let
\[ \omega_p(S,h):=\sup_{s\in[0,h]}\Big(\int_0^{1-s}\norm{S(t+s)-S(t)}^pdt\Big)^{1/p},\quad p\in[1,\infty),\]
with the usual modification for $p=\infty$, denote its $L^p$-modulus of continuity. Then
\[ \abs{S}_{B^\beta_{p,\infty}}:=\sup_{0<h\le 1} h^{-\beta}\omega_p(S,h),\quad 0\le \beta<1,\]
denotes the $B^\beta_{p,\infty}$-seminorm of $S$ and the space $B^\beta_{p,\infty}$ consists of all functions $S\in L^p([0,1];\R^{d\times d})$ with $\abs{S}_{B^\beta_{p,\infty}}<\infty$. For H\"older-continuous $S\in C^\beta(L)$, we have $S\in B^\beta_{p,\infty}$  for any $p$ and $\abs{S}_{B^\beta_{p,\infty}}\le \abs{S}_{B^\beta_{\infty,\infty}}=L$. By the results for Brownian motion \citep{CiesEtal1993} and bounded variation functions \citep{Cohen2003}, any (even discontinuous) It\^o semi-martingale with bounded characteristics is almost surely an element of $B^{1/2}_{p,\infty}$ for any $p\in[1,\infty)$ (but not for $p=\infty$!).

\begin{definition}\label{DefH0H1}
Throughout, we consider the null hypothesis  that the spot covariance has at most rank $r\in\{0,\ldots,d-1\}$. Moreover, we assume a $B^\beta_{p,\infty}$-regularity condition for $\beta\in(0,1)$, $p\in\{2,4\}$, a level $\eps\ge 0$ of idiosyncratic covariance and potentially a spectral gap $\underline\lambda_r> 0$. We set
\begin{align*}
{\mathcal H}_0&:={\mathcal H}_0(r,\beta,L,\eps):=\Big\{(\Sigma_X,\Sigma_Z)\,\Big|\,\sup_{t\in[0,1]}\rank(\Sigma_X(t))\le r,\,\abs{\Sigma_X}_{B^\beta_{2,\infty}}\le L,\,\norm{\Sigma_Z}_{L^2}\le\eps\Big\},\\
{\mathcal H}_0^{\text{gap}}&:={\mathcal H}_0^{\text{gap}}(r,\beta,L,\eps,\underline\lambda_{r}):=\Big\{(\Sigma_X,\Sigma_Z) \,\Big|\,\sup_{t\in[0,1]}\rank(\Sigma_X(t))\le r,\,\inf_{t\in[0,1]}\lambda_r(\Sigma_X(t))\ge\underline\lambda_r,\\
&\qquad\qquad\qquad\qquad\qquad\qquad \,\abs{\Sigma_X}_{B^\beta_{4,\infty}}\le L,\,\norm{\Sigma_Z}_{L^2}\le \eps\Big\}.
\end{align*}

\end{definition}

The role of the parameters is evident and can be seen by the inclusion ${\mathcal H}_0^{gap}(r,\beta,L,\eps,\underline\lambda_{r})\subset {\mathcal H}_0^{gap}(r,\beta',L',\eps',\underline\lambda_r')$ for $\beta\ge\beta'$, $L\le L'$, $\eps\le\eps'$ and $\underline\lambda_r\ge\underline\lambda_r'$.
Sometimes we write $\PP_\Sigma$ for $\Sigma\in{\mathcal H}_0$ to denote the law of model \eqref{EqY} with $\Sigma_X,\Sigma_Z$ from the null hypothesis.


\begin{corollary} \label{CorBias0}
For $\Sigma_X$ with maximal rank $r$ on $I_k$ we have
\begin{align*}
\lambda_{r+1}(\Sigma_X^{kh})&\le \big(\tfrac{2}{\underline\lambda_r}\Delta_2(\Sigma_X,I_k)^2\big)\wedge \Delta_1(\Sigma_X,I_k).
\end{align*}
For the population version of the test statistics $T_{n,h}$ in \eqref{Eqphi0} this gives
\begin{align*}
\sum_{k=0}^{h^{-1}-1}h\lambda_{r+1}(\Sigma_X^{kh})&\le \begin{cases} Lh^\beta,&\text{ under } {\mathcal H}_0(r,\beta,L,\eps),\\ \frac{2}{\underline\lambda_r}L^2h^{2\beta},& \text{ under } {\mathcal H}_0^{gap}(r,\beta,L,\eps,\underline\lambda_r).\end{cases}
\end{align*}
\end{corollary}

\begin{proof}
Apply Proposition \ref{PropPerturb}  with $I=I_k$, $S(t)= \Sigma_X(t)$ to bound $\lambda_{r+1}(\Sigma_X^{kh})$. Without assuming a spectral gap deduce
\begin{align*}
\sum_{k=0}^{h^{-1}-1}h\lambda_{r+1}(\Sigma_X^{kh})&\le  \frac1h\int_0^h\int_0^{1-s} \norm{\Sigma_X(t+s)-\Sigma_X(t)}\,dtds\\
&\le \omega_{L^1}(\Sigma_X,h)\le h^\beta\abs{\Sigma_X}_{B^\beta_{1,\infty}}\le Lh^\beta
\end{align*}
for $\Sigma_X\in {\mathcal H}_0(r,\beta,L,\eps)$ . The spectral gap case follows by the same arguments in terms of $\abs{\Sigma_X}_{B^\beta_{2,\infty}}^2$.
\end{proof}

We are now prepared to derive critical values for the test $\phi_\alpha$ under the null hypotheses without and with a spectral gap. We use the more general setting of observing the process $Y$ and trace back the result for $X$ directly by setting $\eps=0$.

\begin{theorem}\label{ThmNoNoise}
Fix $\alpha\in(0,1)$ and consider the test $\phi_\alpha$ from \eqref{Eqphi0} in terms of $\hat\Sigma_Y^{kh}$. Use the constants $C_{\rho,1},C_{\rho,2}>0$ from Proposition \ref{PropCalcCritVal} below.
\begin{enumerate}
\item Without a minimal spectral gap   take the critical value
\begin{align*}
\kappa_\alpha&= \big(Lh^\beta+\eps\big) \min_{\delta>0}(1+\delta)\Big( 1
  +\big(2C_{d-r,1}+8C_{d-r,2}\big)(nh)^{-1/2}+
 8\delta^{-1}n^{-1/2}\log(\alpha^{-1})\Big).
\end{align*}
Then $\phi_\alpha$ has level $\alpha$ uniformly over ${\mathcal H}_0(r,\beta,L,\eps)$:
$\sup_{\Sigma\in {\mathcal H}_0(r,\beta,L,\eps)}\PP_\Sigma(\phi_\alpha=1)\le \alpha$.

\item Assume the minimal spectral gap $\underline\lambda_r>0$. With critical value
\begin{align*}
\kappa_\alpha
 = \Big(\tfrac{2}{\underline\lambda_r}L^2h^{2\beta}+\eps\Big)\min_{\delta>0}(1+\delta)\Big(& 1+\tfrac{r+4}{2}
  (C_{d-r,1}+4C_{d-r,2}) (nh)^{-1/2}\\
&+ \delta^{-1} (2r+8) \log(\alpha^{-1})n^{-1/2}\Big),
\end{align*}
$\phi_\alpha$ has level $\alpha$ uniformly over ${\mathcal H}_0^{gap}(r,\beta,L,\eps,\underline\lambda_r)$:
$\sup_{\Sigma\in {\mathcal H}_0^{gap}(r,\beta,L,\eps,\underline\lambda_{r})}\PP_\Sigma(\phi_\alpha=1)\le \alpha$.

\end{enumerate}
\end{theorem}

\begin{proof}
This follows from Proposition \ref{PropCalcCritVal} below with $p=2$. Just note $\Delta_p(\Sigma_X,h)\le \abs{\Sigma_X}_{B^\beta_{p,\infty}}h^\beta$  and collect terms.
\end{proof}

The critical values provide uniform and non-asymptotic test levels. They clearly lead to conservative tests because the numerical values are derived from less precise upper bounds. Still, the main structure is clearly discernible: the factors $Lh^\beta+\eps$ or $\frac2{\underline\lambda_r}L^2h^{2\beta}+\eps$ result from the population version in Corollary \ref{CorBias0}, adapted to $\Sigma_Y^{kh}$, while a term of order $(nh)^{-1/2}$ is added to bound the expected values  $\E[\lambda_{r+1}(\hat\Sigma^{kh}_Y-\Sigma^{kh}_Y)]$. The random fluctuations depend on the total sample size and give rise to the factor $n^{-1/2}$. Under a H\"older instead of Besov condition in ${\mathcal H}_0$, we can use $p=\infty$ in Proposition \ref{PropCalcCritVal}. In this case, the factor $(nh)^{-1/2}$ in the term involving $C_{d-r,2}$ can be replaced by $(nh)^{-1}$ and the factor $n^{-1/2}$ in the term involving $\log(\alpha^{-1})$ can be replaced by $n^{-1}$.

We shall see below that the power of the tests is the larger the smaller the critical values are. This relation holds almost independently of the block size $h$.
In view of $n^{-1}\lesssim h\lesssim 1$ we can thus distinguish  two main cases of idiosyncratic noise. For $\eps\lesssim n^{-\beta}$ ($\eps\lesssim \underline\lambda_r^{-1}n^{-2\beta}$ in the spectral gap case), we can choose $h\sim n^{-1}$ and $\kappa_\alpha$ is asymptotically of order $n^{-\beta}$ ($\underline\lambda_r^{-1}n^{-2\beta}$, respectively). Whereas otherwise we can select a larger $h\gg n^{-1}$ so that the term of order $(nh)^{-1/2}$ becomes negligible and $\kappa_\alpha=(1+o(1))\eps$ holds.

The choice of the critical value is based on prior knowledge of the regularity of $\Sigma_X$ and the size of $\Sigma_Z$. Moreover, the limiting law of the test statistics is in general degenerate (deterministic) under ${\mathcal H}_0$. For specific models like the stochastic volatility model in Section \ref{SecStochVol}, data-driven critical values can be established involving quantiles from additional estimation errors.

\subsection{\bf Power and optimal detection rate}

The power analysis shows that the test $\phi_\alpha$ for any choice of critical values is uniformly consistent as $n\to\infty$ over local alternatives whose $(r+1)$st eigenvalues are larger than these critical values. In other words, the $(r+1)$st eigenvalues of $\Sigma_X$ may tend to zero and the tests will still detect them (with asymptotic power one), provided they remain larger than the critical values. This way we establish what is known as a {\it separation rate} between null and alternative in nonparametrics or as a {\it detection rate} in signal processing and learning. Moreover, we establish optimality in the sense that no test can distinguish null hypothesis and local alternatives with $(r+1)$st eigenvalues converging faster to zero.

To define the alternatives formally, we would ideally like to consider covariances $\Sigma_X$ with $\int_0^1\lambda_{r+1}(\Sigma_X(t))dt\ge v_n$ for some detection rate $v_n$. Due to the Riemann-type sum in the tests $\phi_\alpha$, however, we have to use a slightly weaker metric to measure the deviation from zero of the $(r+1)$st eigenvalue. We ask that a box with area $v$ can be placed between the graph of $\lambda_{r+1}(\Sigma_X(t))$ and the $t$-axis, which excludes wild spiky deviations from zero. This may be interpreted as a weak $L^1$-norm over intervals.

\begin{definition}
For $r\in\{0,\ldots,d-1\}$, $\hbar\in (0,1)$, $v>0$ consider the set of alternatives
\begin{align*}
{\mathcal H}_1&:={\mathcal H}_1(r,\hbar,v)
:= \Big\{\Sigma_X\in C([0,1];\R_{spd}^{d\times d})\,\Big|\, \sup_{\abs{I}\ge \hbar}\abs{I}\min_{t\in I}\lambda_{r+1}(\Sigma_X(t))\ge v\Big\}.
\end{align*}
The supremum is taken over all intervals $I\subset[0,1]$ of length at least $\hbar$.
\end{definition}

Remark that the alternative ${\mathcal H}_1$  neither involves a $\beta$-smoothness constraint nor a spectral gap condition.

\begin{theorem}\label{ThmPower}
Consider for sample size $n$ the tests $\phi_{\alpha,n}$ from \eqref{Eqphi0}  in terms of $\hat\Sigma_Y^{kh}$ with block sizes $h_n$ and some critical values $\kappa_{\alpha,n}$. Assume $h_n\to 0$ as $n\to\infty$ and that the number of observations per block satisfies $nh_n\ge 2(r+1)C_{d,r+1}$ with the constant $C_{d,\ell}$ from Corollary \ref{CorEigenvalueTriangular} below. Then $\phi_{\alpha,n}$ is asymptotically consistent  over the local alternatives ${\mathcal H}_1(r,\hbar_n,v_n)$ provided $\hbar_n/h_n\to\infty$ and the rate $v_n$ is larger than a constant multiple of $\kappa_{\alpha,n}$:
\[ \exists C>0:\,\lim_{n\to\infty} \inf_{\Sigma_X\in {\mathcal H}_1(r,\hbar_n,C\kappa_{\alpha,n})}\PP_{\Sigma_X}(\phi_{\alpha,n}=1)=1.\]
\end{theorem}

\begin{remark}
The idiosyncratic process $Z$ may be zero  or non-zero  under the alternative because of $\Sigma_Y(t)=\Sigma_X(t)+\Sigma_Z(t)\ge\Sigma_X(t)$. In case $\eps'=\inf_t\lambda_{min}(\Sigma_Z(t))>0$, the consistency result can be shown to hold already for ${\mathcal H}_1(r,\hbar_n,C(\kappa_{\alpha,n}-\eps'))$. In case $nh_n\to\infty$ it suffices to require $\hbar_n\gtrsim h_n$ only. The imposed lower bound on the number $nh_n$ of observations per block is finite, but quite pessimistic (it grows exponentially in $d$). We were not able to establish the result under the minimal possible bound $nh_n\ge r+1$. Simulation results indicate in any case that a choice $nh_n=C(r+1)$ with $C$ of the order 10 produces  stable results, see also Example \ref{ExRank2} below.
\end{remark}

\begin{proof}
Consider $A_{jk}=n\int_{I_{jk}}\Sigma_Y(t)dt$ with $I_{jk}=[kh_n+(j-1)/n,kh_n+j/n]$ for $j=1,\ldots,nh_n$ and $k=0,\ldots,h_n^{-1}-1$. For $\Sigma_X \in{\mathcal H}_1(r,\hbar_n,v_n)$ choose $a,b\in[0,1]$  with $\lambda_{r+1}(\Sigma_X(t))\ge (b-a)^{-1}v_n$ for $t\in[a,b]$ and $b-a\ge \hbar_n$. Then by $\Sigma_Y(t)\ge\Sigma_X(t)$ and Lemma \ref{Lemtracer} there are at least $K_n':=\floor{(b-a)/h_n}\ge 1$ blocks $I_k$ where
\[ \lambda_{r+1}(A_{jk})\ge (d-r)^{-1}n\int_{I_{jk}} \trace_{>r}(\Sigma_X(t))\,dt\ge (d-r)^{-1}(b-a)^{-1}v_n\]
holds for all $j$. Using this bound in Corollary \ref{CorEigenvalueTriangular} below, we obtain with $J_0=\ceil{nh_n/C_{d,r+1}}\ge 2(r+1)$ (by assumption) for any $\tau>0$
\begin{align*}
\PP_{\Sigma_X}\Big(\sum_{k=0}^{h_n^{-1}-1} \lambda_{r+1}\Big(\hat\Sigma_Y^{kh_n}\Big)\le \tau K_n'(d-r)^{-1}(b-a)^{-1}v_n\Big)\le
\Big(\bar C_{d,r+1}\tau\Big)^{K_n'(J_0-r)/2}.
\end{align*}
We take $\tau=(2\bar C_{d,r+1})^{-1}$ and insert
$h_nK_n'(b-a)^{-1}\ge 1-h_n/\hbar_n$
to arrive at
\begin{align*}
\PP_{\Sigma_X}\Big(\sum_{k=0}^{h_n^{-1}-1} h_n\lambda_{r+1}\Big(\hat\Sigma_Y^{kh_n}\Big)\le (2\bar C_{d,r+1})^{-1} (d-r)^{-1}(1-h_n/\hbar_n)v_n\Big)\le
2^{-K_n'(J_0-r)/2}.
\end{align*}
Because of $\hbar_n/h_n\to\infty$ the number $K_n'$ of blocks tends to infinity and hence
\begin{align*}
\lim_{n\to\infty}\inf_{\Sigma_X\in{\mathcal H}_1(r,h_n',v_n)}\PP_{\Sigma_X}\Big(\sum_{k=0}^{h_n^{-1}-1} h_n\lambda_{r+1}\Big(\hat\Sigma_Y^{kh_n}\Big)> cv_n\Big)=1
\end{align*}
follows for any $c\in(0,(2\bar C_{d,r+1})^{-1}(d-r)^{-1})$.  This is the assertion when setting $C=c^{-1}$, $v_n=C\kappa_{\alpha,n}$.
\end{proof}

The smaller $h_n$, the smaller the order of the critical values in Theorem \ref{ThmNoNoise}. The bound under the alternative merely requires $h_n\ge Cn^{-1}$ for some sufficiently large constant $C>0$. In conclusion we can choose $h_n\thicksim n^{-1}$ minimal without losing asymptotic power. Remark in this case that even for constant $\Sigma_X$ the test statistic $\sum_k h_n\lambda_{r+1}(\hat\Sigma_X^{kh_n})$ does not estimate $\lambda_{r+1}(\Sigma_X)$, but rather the expected $(r+1)$st eigenvalue of the corresponding Wishart distribution. This expected value grows in $\lambda_{r+1}(\Sigma_X)$ and thus suffices for inference. Combining Theorems \ref{ThmNoNoise} and \ref{ThmPower}, the tests $\phi_{\alpha,n}$ with $h_n\thicksim n^{-1}$ thus establish the following signal detection rate.

\begin{corollary}
Given $n$ observations, there are  tests $\phi_{\alpha,n}$ of ${\mathcal H}_0(r,\beta,L,\eps_n)$ or ${\mathcal H}_0^{\text{gap}}(r,\beta,L,\eps_n,\underline\lambda_{r,n})$, respectively, versus ${\mathcal H}_1(r,\hbar_n,v_n)$ that have uniform level $\alpha$  and are uniformly consistent over the alternatives for $n\to\infty$ if $n\hbar_n\to\infty$ and for some suitably large constant $C>0$
\[ v_n=C\Big(\big(Ln^{-\beta}+\eps_n\big)\wedge \big(\underline\lambda_{r,n}^{-1}L^2n^{-2\beta}+\eps_n\big) \Big).\]
The parameters $\eps_n$ and $\underline\lambda_{r,n}$  may vary arbitrarily with $n$.
\end{corollary}

\begin{remark}\label{RemDrift}
If the pure Brownian martingale model \eqref{EqX} is generalised to the semi-martingale model
\[ dX(t)=b_X(t)dt+\sigma_X(t)dB(t),\]
including a bounded adapted drift $b_X(t)\in\R^d$, then the increments of $X$ involve an additional term of order ${\mathcal O}_P(n^{-1})$ such that $\hat\Sigma_X^{kh}$ and the corresponding test statistics $T_{n,h}$ in \eqref{Eqphi0} are perturbed by a bias term of order ${\mathcal O}_P(n^{-1})$. In this case we should employ critical values larger than $n^{-1}$, that is choose $h$ appropriately, to remain robust against these perturbations or alternatively use a mean-corrected version of the blockwise covariance estimator $\hat\Sigma_X^{kh}$.

In case $\sigma_X(t)$ is stochastic and dependent on the driving Brownian motion $B$ in an adapted way, a standard procedure would be to approximate $\sigma_X(t)$ by $\sigma_X(kh)$ on $I_k$ and  to work conditionally on the underlying filtration ${\mathcal F}_{kh}$, where the covariance can then be treated as deterministic and constant. At a population level we would work with $\Sigma_X(kh)$ instead of $\Sigma_X^{kh}$. We expect that pursuing this idea we can at best achieve similar results to those without a spectral gap. A direct analysis of $\lambda_{r+1}(\hat\Sigma_X^{kh})$ might yield better results, but would certainly require to develop different techniques.
\end{remark}

That the above rates cannot be improved follows basically from Example \ref{ExRank}. In Appendix \ref{AppTechnical3}  the following non-asymptotic lower bound or rather non-identifiability result is established  for the case $\eps=0$.

\begin{proposition}\label{PropLB}
Consider the model of observing $X$ in \eqref{EqX} at times $i/n$, $i=0,\ldots,n$, and arbitrary parameters $1\le r<d$, $\beta\in(0,1)$, $L>0$, $\underline\lambda_r> 0$ as well as $\hbar\in(0,1]$.
Then there is no non-trivial test between ${\mathcal H}_0$ and ${\mathcal H}_1$:
\begin{align*} \forall\text{ tests }\phi:&\, \sup_{\Sigma\in {\mathcal H}_0(r,\beta,L,0)}\E_\Sigma[\phi]+ \sup_{\Sigma\in {\mathcal H}_1(r,\hbar,cLn^{-\beta})}\E_\Sigma[1-\phi]=1,\\
\forall\text{ tests }\phi:&\, \sup_{\Sigma\in {\mathcal H}_0^{gap}(r,\beta,L,0,\underline\lambda_{r})}\E_\Sigma[\phi]+ \sup_{\Sigma\in {\mathcal H}_1(r,\hbar,(\underline\lambda_r^{-1}c^2L^2n^{-2\beta})\wedge(cLn^{-\beta}))}\E_\Sigma[1-\phi]=1
\end{align*}
whenever $c\le2^{-5/2}\pi^{-1}$.
\end{proposition}

We conclude that the optimal detection rate is indeed
\[ v_n=(\underline\lambda_r^{-1}L^2n^{-2\beta})\wedge(Ln^{-\beta})\]
whenever $\eps$ is smaller than $v_n$. It is quite remarkable that this rate is so much faster than the estimation rate $n^{-\beta/(2\beta+1)}$ for the spot covariance matrix $\Sigma_X(t)$. A fundamental reason is, of course, the heteroskedasticity in the estimation error which becomes small in directions where $\Sigma_X(t)$ is small. Compared to the $n^{-1/2}$-estimation rate for positive realised integrated eigenvalues, this provides a quantitative version of  the super-efficiency around zero eigenvalues \citep[Remark 3]{ait2019}. It should be remarked that the case $r=0$ has been excluded because it leads to the completely degenerate model $\Sigma_X(t)=0$ and thus $X(t)=0$ under ${\mathcal H}_0$, which can be tested perfectly.

\begin{figure}[t]
\centering
\includegraphics[width=0.48\textwidth]{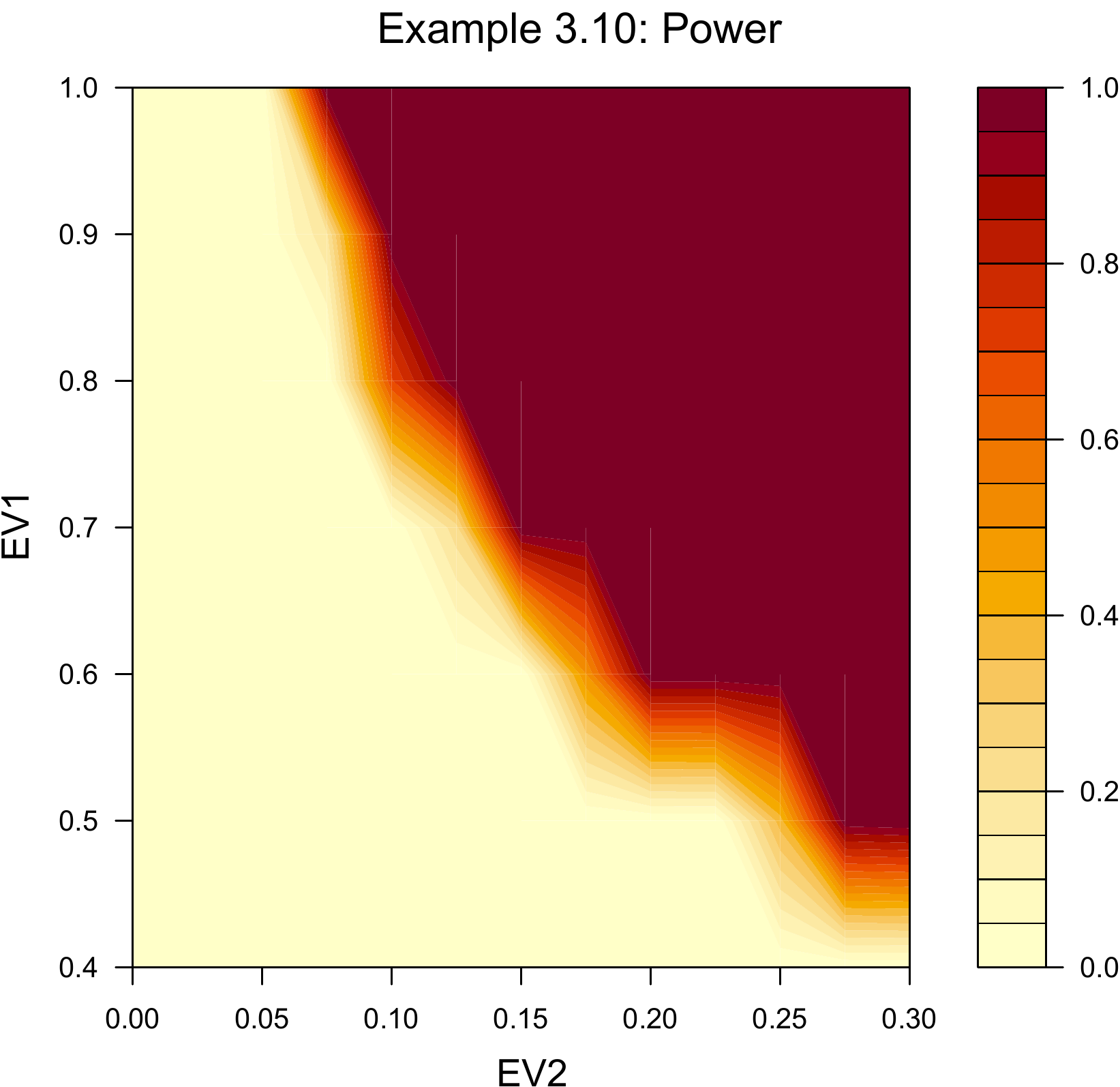} \includegraphics[width=0.48\textwidth]{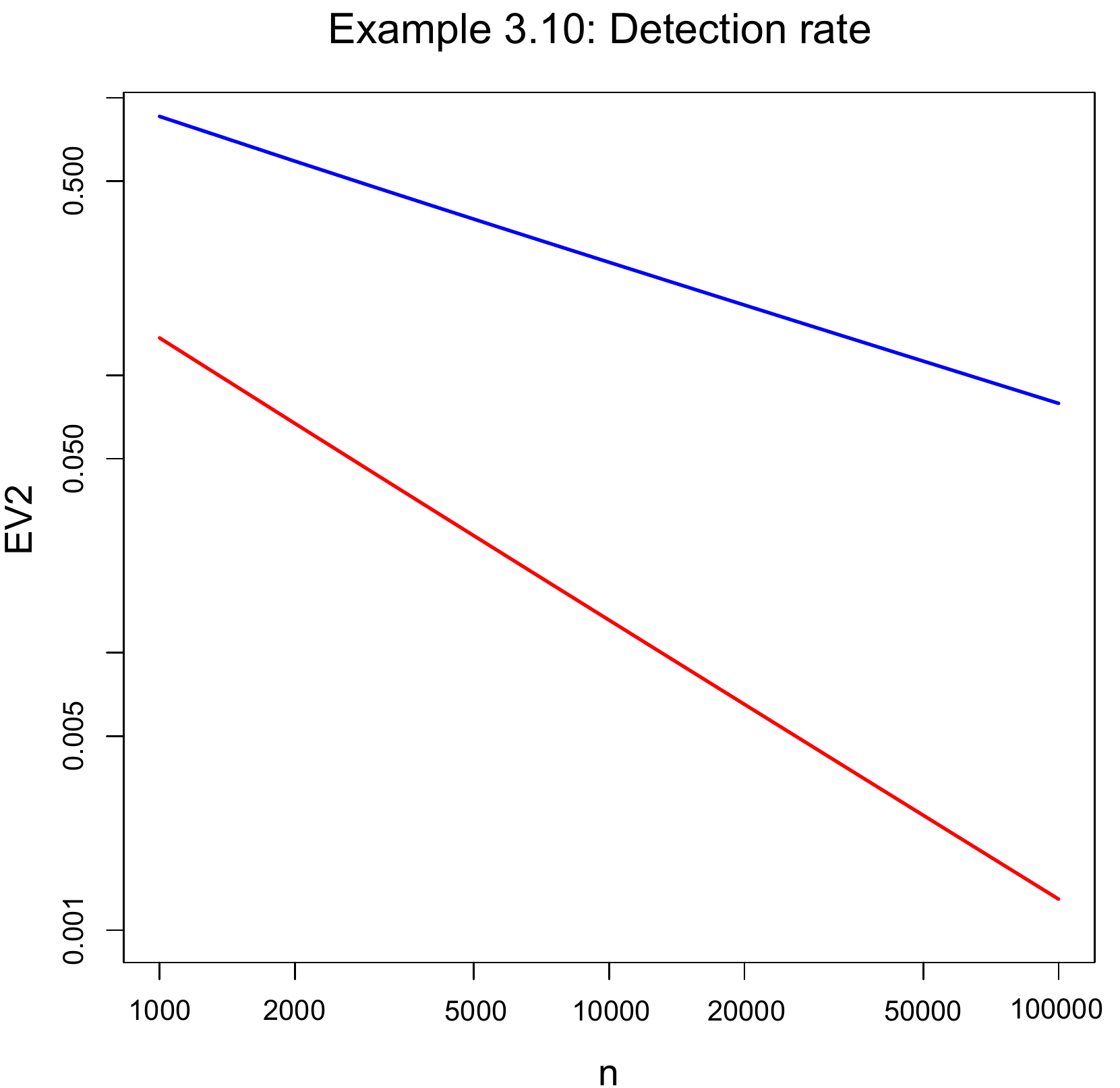}

 \caption{Left: Power of the rank test as a function of signal strength (EV2) and spectral gap (EV1) in Example \ref{ExRank2}. Right: 50\% detection rates (log-log plot) with spectral gap (red) and without (blue).}\label{FigPower}
\end{figure}

\begin{example}\label{ExRank2}
We generalise Example \ref{ExRank} with $\beta=1/2$ to illustrate  the impact of the spectral gap and the magnitude of the second eigenvalue of the spot covariance matrix on the power and on the detection rate of the test. With $v_1(t)=(\underline\lambda_1^{1/2},(h/\underline\lambda_1)^{1/2}\sin(2\pi t/h))^\top$ and $v_2(t)=((h/\underline\lambda_1)^{1/2}\sin(2\pi t/h), -\underline\lambda_1^{1/2})^\top$, we set
$\Sigma_X(t)=v_1(t)v_1(t)^\top+\gamma v_2(t)v_2(t)^\top$.
By orthogonality, $v_1(t)$ and $v_2(t)$ are both  eigenvectors of $\Sigma_X(t)$ with eigenvalues $\lambda_1(\Sigma_X(t))=\underline\lambda_1+(h/\underline\lambda_1)\sin^2(2\pi t/h)$ and $\lambda_2(\Sigma_X(t))=\gamma\lambda_1(\Sigma_X(t))$ for $\gamma\in[0,1]$. By varying $\gamma$, we generate spot covariances under the alternative ${\mathcal H}_1: r=2$ with different signal strengths $h^{-1}\int_0^h\lambda_2(\Sigma_X(t))\,dt=\gamma(\underline\lambda_1+h/(2\underline\lambda_1))$.

The contour plot in Figure \ref{FigPower}(left) displays the power, that is the probability that $\phi_\alpha$ with $\alpha=0.1$ from \eqref{Eqphi0} rejects. We use the minimum of the critical values $\kappa_\alpha$ from Theorem \ref{ThmNoNoise}, depending on the spectral gap $\underline\lambda_1$. Warmer colours indicate a higher power, which is shown as a function of the average second eigenvalue (signal strength) on the $x$-axis and the spectral gap $\underline\lambda_1$ on the $y$-axis. We use $n=2000$ observations in $1000$ Monte Carlo iterations with $h=0.02$. The power increases with the signal (average second eigenvalue) and with the spectral gap. The phase transition from almost zero to perfect power is very fast, especially for larger spectral gaps, in line with the mathematical analysis.

Figure \ref{FigPower}(right) shows for the same example a log-log-plot of two signal detection rates as a function of the sample size $n$. For the case of a spectral gap ($\underline\lambda_1=1$, red) and without a spectral gap (blue) the value EV2 indicates the  average second eigenvalue at which the test accepts and rejects with equal probability $50\%$. This corresponds to the standard classification boundary in learning. Given the sharp phase transitions, compare Figure \ref{FigPower}(left), the values of EV2 are almost perfectly given by the respective critical values in Theorem \ref{ThmNoNoise}. Consequently, the power laws in the detection rate $v_n$ meet well the finite sample classification boundaries. We have chosen $nh=40$ observations per block for each $n$ and calculated the power in $1000$ Monte Carlo simulations for each case.
\end{example}

\begin{remark}
Eigenvalue detection is intrinsically different from standard signal detection problems. When testing whether the regression function $f$ in nonparametric regression is zero versus local alternatives with $\norm{f}_{L^p}\ge v_n$, a smoothness or sparsity condition needs to be imposed in the alternative \citep{ingster2012}. In our case, the smoothness condition appears under the null hypothesis $\lambda_{r+1}(\Sigma_X(\cdot))=0$, but not under the alternative. Moreover, we face a one-sided testing problem because always $\lambda_{r+1}(\Sigma_X(t))\ge 0$ which is basically trivial in nonparametric regression. A major difference is, of course, that our null hypothesis $\lambda_{r+1}(\Sigma_X(\cdot))=0$ is composite, allowing for a submanifold of covariance matrix functions $\Sigma_X$.

There seems to be, however, also a deeper reason for the inversion of roles of ${\mathcal H}_0$ and ${\mathcal H}_1$, which is the concavity of the $\trace_{>r}$-functional when expressing ${\mathcal H}_0$ as $\trace_{>r}(\Sigma_X(\cdot))=0$ and ${\mathcal H}_1$ as $\trace_{>r}(\Sigma_X(t))>0$ (in a weak $L^1$-norm over $t$). The concavity is inherent in the explanations for the examples in Section \ref{SecSetting} as well as in several proofs. This should be
compared with the convexity of the $L^p$-norms and standard shape constraints in regression problems, compare  the discussion in \citet{juditsky2002}. From a scientific inference perspective the null hypothesis should always consider the simpler and less complex setting ({\it Occam's razor principle}), which clearly justifies our approach even though the information geometry may favour the opposite.
\end{remark}

\subsection{\bf A rank estimator}

The rank test translates to an estimator $\hat r$ of the rank $r\in\{0,\ldots,d\}$. 
We denote by $\phi_{\alpha}^{(r)}$ and $\kappa_{\alpha}^{(r)}$  the level-$\alpha$ test $\phi_{\alpha}$ for $H_0:\rank(\Sigma(t))\le r$ from \eqref{Eqphi0} with critical value $\kappa_{\alpha}$. A natural rank estimator uses the minimal rank bound that is accepted by sequential testing:
\begin{equation}\label{EqDefhatr} \hat r:=\inf\{j=0,\ldots,d-1\,|\,\phi_{\alpha_j}^{(j)}=0\}\wedge d
\end{equation}
with suitable levels $\alpha_j$. Setting $\hat\lambda_j:=\sum_{k=0}^{h^{-1}-1}h\lambda_j(\hat\Sigma_Y^{kh})$, this estimator can be rewritten in form of the optimisation problem
\[ \hat r=\argmin_{0\le r\le d}\Big(\sum_{j=r+1}^d \hat\lambda_j+\sum_{j=0}^{r-1} \kappa_{\alpha_j}^{(j)}\Big),\]
provided the critical values $\kappa_{\alpha_j}^{(j)}$ are non-decreasing in $j$. Indeed, since subtracting the trace  $\hat\lambda_1+\cdots+\hat\lambda_d$ does not change the minimiser, we have
\begin{equation}\label{EqRankCrit} \argmin_{0\le r\le d}\Big(\sum_{j=r+1}^d \hat\lambda_j+\sum_{j=0}^{r-1} \kappa_{\alpha_j}^{(j)}\Big)=\argmin_{0\le r\le d}\sum_{j=0}^{r-1} \big(\kappa_{\alpha_j}^{(j)}-\hat\lambda_{j+1}\big).
\end{equation}
The summand on the right is almost surely strictly increasing in $j$ and must thus be  negative for $j\le \hat r-1$ and be positive for $j\ge \hat r$. Hence,
\[ \hat r=\inf\big\{j=0,\ldots,d-1\,\big|\,\kappa_{\alpha_j}^{(j)}-\hat\lambda_{j+1}\ge 0 \big\}\wedge d=\inf\{j=0,\ldots,d-1\,|\,\phi_{\alpha_j}^{(j)}=0\}\wedge d
\]
follows.

Let us specify this for the tests $\phi_\alpha^{(r)}={\bf 1}(\hat\lambda_{r+1}>\kappa_{\alpha_r}^{(r)})$ from Theorem \ref{ThmNoNoise} without a spectral gap assumption. We fix a level $\alpha\in(0,1)$ and set
\[ \kappa_\alpha:=\kappa_\alpha^{(0)}=\big(Lh^\beta+\eps\big)\min_{\delta>0}(1+\delta)\Big( 1
  +\big(2C_{d,1}+8C_{d,2}\big)(nh)^{-1/2}+
 8\delta^{-1}n^{-1/2}\log(\alpha^{-1})\Big).
\]
Then $\kappa_\alpha^{(r)}\le \kappa_\alpha$ holds for all $r=0,\ldots,d-1$, but the difference is only in the second order term via the remaining dimension $d-r$. This leads to the choice
\begin{equation}\label{Eqhatr}
\hat r=\argmin_{0\le r\le d}\Big(r\kappa_\alpha+\sum_{j=r+1}^d \hat\lambda_j\Big).
\end{equation}

\begin{proposition}\label{PropRankDetect}
The rank estimator $\hat r$ from \eqref{Eqhatr} satisfies for $\beta\in(0,1)$, $L>0$, $\eps>0$ non-asymptotically
\[ \max_{r=0,\ldots,d-1}\sup_{\Sigma\in {\mathcal H}_0(r,\beta,L,\eps)}\PP_\Sigma(\hat r>r)\le \alpha.\]
Moreover, asymptotically for $h_n,\hbar_n\to 0$ as $n\to\infty$, $\hbar_n/h_n\to\infty$ and $nh_n\ge 2(r+1)C_{d,r+1}$, with $C_{d,r+1}$ from Corollary \ref{CorEigenvalueTriangular}, there is a constant $C>0$ such that
\[ \lim_{n\to\infty}\max_{r=1,\ldots,d}\sup_{\Sigma_X\in {\mathcal H}_1(r,\hbar_n,C(Lh_n^\beta+\eps_n))}\PP_\Sigma(\hat r<r)=0,\]
where $\alpha=\alpha_n$ may tend to zero, provided $\alpha_n\ge \exp(-cn^{1/2})$ holds for some constant $c>0$.
\end{proposition}

\begin{proof}
From Theorem \ref{ThmNoNoise} we obtain directly
\[ \max_{r=0,\ldots,d-1}\sup_{\Sigma\in {\mathcal H}_0(r,\beta,L,\eps)}\PP_\Sigma(\hat r>r)\le \max_{r=0,\ldots,d-1}\sup_{\Sigma\in {\mathcal H}_0(r,\beta,L,\eps)}\PP_\Sigma(\hat \lambda_{r+1}> \kappa_\alpha)\le \alpha
\]
because of  Representation \eqref{EqRankCrit} and $\kappa_\alpha\ge \kappa_\alpha^{(r)}$.

In the asymptotic setting, we have $\kappa_\alpha\lesssim (Lh_n^\beta+\eps_n)$ because of $nh_n\ge 1$, $\log(\alpha^{-1})\lesssim n^{-1/2}$. Hence, the definition of $\hat r$ and Theorem \ref{ThmPower} gives for some sufficiently large constant $C>0$ that
\[ \max_{r=1,\ldots,d}\sup_{\Sigma_X\in {\mathcal H}_1(r,\hbar_n,C(Lh_n^\beta+\eps_n))}\PP_{\Sigma_X}(\hat r<r) \le \max_{r=1,\ldots,d}\sup_{\Sigma_X\in {\mathcal H}_1(r,\hbar_n,C(Lh_n^\beta+\eps_n))} \PP_{\Sigma_X}(\hat\lambda_{r}\le\kappa_\alpha)
\]
tends to zero.
\end{proof}

Our rank estimator $\hat r$ strictly controls the probability of selecting a too large rank. If we assume $\eps_n\to 0$, then the probability of choosing a too small rank for any fixed model tends to zero asymptotically because any continuous $\Sigma_X$ with $\int_0^1\lambda_{r+1}(\Sigma_X(t))dt>0$ lies in some ${\mathcal H}_1(r,\delta_1,\delta_2)$ for sufficiently small $\delta_1,\delta_2>0$. The uniform control of Proposition \ref{PropRankDetect} is much stronger than standard pointwise consistency results in the literature. For $\alpha_n\thicksim e^{-cn^{1/2}}$ with sufficiently small $c>0$ we conjecture that even the joint error probability decays as fast as $\PP(\hat r\not=r)\le e^{-cn^{1/2}}$. Similar results can be obtained in the spectral gap case with the asymptotics $\underline\lambda_r\to 0$.

Compared to the \citet{BaiNg2002} and \citet{ait2017} estimators, designed in a slightly different setting, we observe that an exact form of the penalty is provided by bounding the error of overestimating the rank by $\alpha$ non-asymptotically. We believe that the approach \eqref{EqDefhatr} of determining $\hat r$ by sequential testing is both, from a conceptual and a practical point of view, very attractive, while the optimisation formulation \eqref{Eqhatr} allows a better comparison with other existing rank detection methods.

\section{Stochastic volatility models} \label{SecStochVol}

For certain stochastic volatility models we can access the bias bounds $Lh^\beta$ and $L^2h^{2\beta}$ for our test statistics by an estimation procedure.
Assume that the covariance can be modeled by a continuous semi-martingale
\begin{equation}\label{EqStochVol}
d\Sigma_X(t)=b(t)dt+\Gamma(t)dB'(t),\quad t\in[0,1],\qquad \Sigma(0)=\Sigma_0,
\end{equation}
satisfying the following assumptions, where ${\mathcal L}(V,W)$ denotes all linear maps between vector spaces $V$ and $W$.

\begin{assumption}\label{AssStochVol}
Fix $p\in\{1,2\}$. $B'$ is a $d'$-dimensional Brownian motion, $b(t)\in\R^{d\times d}_{sym}$ is an $L^{2p}$-bounded, adapted covariance drift and  $\Gamma(t)\in {\mathcal L}(\R^{d'},\R^{d\times d}_{sym})$ is an $L^{4p}$-bounded, adapted covariance diffusivity. To avoid clumsy vectorisations, we interpret $\Gamma(t)dB'(t)$ as a matrix-valued integrator with linear combinations of the coordinates $dB'_i(t)$ in each entry.

The regularity conditions $\E[\norm{\Gamma(t)-\Gamma(s)}^{2p}]^{1/2p}\le L_\Gamma\abs{t-s}^{\beta_\Gamma}$, $\E[\norm{b(t)-b(s)}^{2p}]^{1/2p}\le L_b\abs{t-s}^{\beta_b}$ hold for all $t,s$ and some constants $L_\Gamma,L_b>0$, $\beta_\Gamma>1/2$, $\beta_b>0$. All processes $b$, $\Gamma$ and $B'$ as well as the ${\mathcal F}_0$-measurable initial condition $\Sigma_0\in\R^{d\times d}_{spd}$ in $L^{2p}$ are  defined on the same filtered probability space as $B$ and $X$. Moreover,  $B'$ and $B$  are independent. Model \eqref{EqStochVol} enforces the symmetric matrix $\Sigma_X(t)$ to be positive semidefinite for all $t$.
\end{assumption}


In this section we consider only the case where we observe $X$ directly, i.e. $Y=X$ and $\eps=0$. In the case  $\eps>0$ the subsequent approach will generally overestimate the $p$-variations and  lead to more conservative tests.

The first limit result shows that the bias bound from Corollary \ref{CorBias0} can be expressed asymptotically by the {\it normed p-variation}
\[NV^{(p)}:=\int_0^1\rho_p(\Gamma(t))dt\text{, where }\rho_p(\Gamma(t)):=\E_Z[\norm{\Gamma(t) Z}^p]\]
and the expectation is taken with respect to an independent random vector $Z\sim N(0,I_{d'})$ (in general, $NV^{(p)}$ is thus random via $\Gamma(t)$).

\begin{proposition}\label{PropNQV}
For a semi-martingale $\Sigma_X$ satisfying \eqref{EqStochVol} under Assumption \ref{AssStochVol} and for $p\in\{1,2\}$ we  have
\[ \Delta_p(\Sigma_X,h)^p:=\sum_{k=0}^{h^{-1}-1}h\Delta_p(\Sigma_X,I_k)^p=h^{p/2}\Big( \tfrac {8}{(p+2)(p+4)}NV^{(p)}+{\mathcal O}_{L^2}(h^{1/2})\Big).
\]
\end{proposition}

\begin{proof}
Bounding the drift and using the regularity of $\Gamma$ we derive
\begin{align*}
&\int_{I_k}\int_{I_k}\norm{\Sigma_X(t)-\Sigma_X(s)}^pdsdt = 2\int_{kh}^{(k+1)h}\int_{kh}^t\bnorm{\int_s^t\Gamma(u)dB'(u)+{\mathcal O}_{L^{2p}}(h)}^pdsdt\\
&=2\int_{kh}^{(k+1)h}\int_{kh}^t\bnorm{\Gamma(kh)(B'(t)-B'(s))+{\mathcal O}_{L^{2p}}(h^{\beta_\Gamma+1/2}+h)}^2dsdt.
\end{align*}
Hence, taking conditional expectation and arguing for $p=2$ via Cauchy-Schwarz inequality with $\norm{B'(t)-B'(s)}={\mathcal O}_{L^4}(h^{1/2})$, $\beta_{\Gamma}+1/2\ge 1$ we find
\begin{align*}
&\sum_{k=0}^{h^{-1}-1}\E\Big[\frac1{h}\int_{I_k}\int_{I_k}\norm{\Sigma_X(t)-\Sigma_X(s)}^pdsdt\,\Big|\,{\mathcal F}_{kh}\Big]\\
&=\sum_{k=0}^{h^{-1}-1}\rho_p(\Gamma(kh))\frac2h\int_0^h\int_0^t(t-s)^{p/2}dsdt+ {\mathcal O}_{L^2}(h^{(p+1)/2})\\
&=h^{p/2}\Big(\frac{8}{(p+2)(p+4)}NV^{(p)}+{\mathcal O}_{L^2}(h^{\beta_\Gamma}+h^{1/2})\Big),
\end{align*}
because the mapping $\Gamma\mapsto\rho_p(\Gamma)^{1/p}$ is Lipschitz continuous so that with a uniform bound on $\rho_p(\Gamma(t))$ we deduce $\abs{\rho_p(\Gamma(t))-\rho_p(\Gamma(s))}={\mathcal O}_{L^2}(\abs{t-s}^{\beta_\Gamma})$. Moreover, we directly get uniformly in $k$
\[ \E\Big[\Big(\frac1{h}\int_{I_k}\int_{I_k}\norm{\Sigma_X(t)-\Sigma_X(s)}^pdsdt\Big)^2\Big]={\mathcal O}(h^{p+2}).\]
We conclude by profiting from the martingale difference structure
\begin{align*}
&\E\Big[\Big(\frac1{h}\sum_{k=0}^{h^{-1}-1} \Big(\int_{I_k}\int_{I_k}\norm{\Sigma_X(t)-\Sigma_X(s)}^p dsdt- \E\Big[\int_{I_k}\int_{I_k}\norm{\Sigma_X(t)-\Sigma_X(s)}^pdsdt\,\Big|\,{\mathcal F}_{kh}\Big]\Big)\Big)^2\Big]\\
&\le \sum_{k=0}^{h^{-1}-1}\E\Big[\Big(\frac1{h}\int_{I_k}\int_{I_k}\norm{\Sigma_X(t)-\Sigma_X(s)}^pdsdt\Big)^2\Big]
={\mathcal O}(h^{p+1}).
\end{align*}
Because of $\beta_\gamma\ge 1/2$ the total $L^2$-approximation error is of order ${\mathcal O}(h^{(p+1)/2})$.
\end{proof}

We shall estimate $NV^{(p)}$ on a coarser grid with block length $h'$. For scalar models it is known that estimators for the volatility of volatility have at best rate $n^{-1/4}$ \citep{Hoffmann2002}. If no subtle bias corrections are employed, the standard rate is (almost) $n^{-1/5}$, see \citet{Vetter2015} and the discussion therein. Since we cannot ensure differentiability of the matrix norm without further assumptions, we are content with a consistent estimator achieving at best the rate $n^{-1/6}$.  To render the influence of the finite variation part in \eqref{EqStochVol} sufficiently small, we use  second differences. Techniques from discretisation of processes give the following central limit theorem, proved in Appendix \ref{AppStochVol}.

\begin{theorem}\label{ThmNQVCLT}
The normed $p$-variation estimator
\begin{equation}\label{EqDefNpVhat}
 \widehat{NV}^{(p)}_{h',n}:=\frac{3h'}{(2h')^{p/2}}\sum_{k=0}^{(3h')^{-1}-1} \norm{\hat\Sigma_X^{(3k+2)h'}-2\hat\Sigma_X^{(3k+1)h'}+\hat\Sigma_X^{3kh'}}^p
 \end{equation}
satisfies in Model \eqref{EqStochVol} under Assumption \ref{AssStochVol} for block sizes $h'\to 0$ with $h'n^{1/3}\to\infty$
\[ \Big(3h'\int_0^1 \big(\rho_{2p}(\Gamma(t))-\rho_p(\Gamma(t))^2\big)\,dt\Big)^{-1/2}\Big(\widehat{NV}^{(p)}_{h',n}-NV^{(p)}\Big) \xrightarrow{d} N(0,1).
\]
\end{theorem}

To obtain a feasible central limit theorem we need to estimate the variance term $\int_0^1 (\rho_{2p}(\Gamma(t))-\rho_p(\Gamma(t))^2)dt$ consistently. This can be accomplished using the $2p$-powers of norms as well as the product of $p$-powers on adjacent blocks.

\begin{proposition}\label{PropVarhat}
Grant Assumption \ref{AssStochVol} for $p\in\{1,2\}$ as well as $h'\to 0$ with $h'n^{1/2}\to\infty$. Then $\widehat{NV}_{h',n}^{(2p)}$ from \eqref{EqDefNpVhat}  is consistent:
\[ \widehat{NV}_{h',n}^{(2p)} \xrightarrow{\PP} \int_0^1\rho_{2p}(\Gamma(t))\,dt.
\]
Moreover, the bipower normed $p$-variation estimator $\widehat{BNV}^{(p)}_{h',n}$ given by
\begin{align*}
\frac{6h'}{(2h')^p}&\sum_{k=0}^{(6h')^{-1}-1} \norm{\hat\Sigma_X^{(6k+2)h'}-2\hat\Sigma_X^{(6k+1)h'}+\hat\Sigma_X^{6kh'}}^p \norm{\hat\Sigma_X^{(6k+5)h'}-2\hat\Sigma_X^{(6k+4)h'}+\hat\Sigma_X^{(6k+3)h'}}^p
\end{align*}
estimates the squared normed $p$-variation consistently:
\begin{align*}
\widehat{BNV}^{(p)}_{h',n}
& \xrightarrow{\PP} \int_0^1\rho_p(\Gamma(t))^2dt.
\end{align*}
\end{proposition}

The last two results allow to replace the term $(Lh^\beta)^p$, $p\in\{1,2\}$, in the critical values of Theorem \ref{ThmNoNoise} by $\tfrac {8}{(p+2)(p+4)}\widehat{NV}^{(p)}_{h',n}h^{p/2}$ plus a Gaussian quantile depending on the estimated variance of $\widehat{NV}^{(p)}_{h',n}$. For $h'\gg (nh)^{-1}$ the lower order terms in the critical values of Theorem \ref{ThmNoNoise} become negligible and we obtain level-$\alpha$ tests with quite simple data-driven critical values.

\begin{corollary}\label{Corkappahat1}
Assume Model \eqref{EqStochVol} under Assumption \ref{AssStochVol} for $p=1$. Consider block sizes $h,h'\to 0$ with $h=o(h')$, $h'n^{1/3}\to\infty$ and $nhh'\to\infty$.
Then the test $\phi_\alpha=\phi_{\alpha,n,h,h'}$ from \eqref{Eqphi0} with fixed $\alpha\in(0,1)$ and critical value
\[ \kappa_\alpha=\kappa_{\alpha,n,h,h'}=\tfrac {8}{15}h^{1/2}\Big(\widehat{NV}_{h',n}^{(1)}+ \Big(3h'\big(\widehat{NV}_{h',n}^{(2)}-\widehat{BNV}_{h',n}^{(1)}\big)\Big)^{1/2}q_{1-\alpha;N(0,1)}\Big)\]
has asymptotic level $\alpha$ on the (possibly random) null hypothesis ${\mathcal H}_0=\{\sup_{t\in[0,1]}\rank(\Sigma_X(t))\le r\}$:
\[ \limsup_{n\to\infty}\PP\Big(\{\phi_{\alpha,n,h,h'}=1\}\cap{\mathcal H}_0\Big)\le\alpha.\]
\end{corollary}

\begin{proof}
By Proposition \ref{PropCalcCritVal} below (set $\Sigma_Z=0$, $p=2$, $\delta\to 0$ with $\delta^{-1}\le h^{-1/2}$), we have on the event ${\mathcal H_0}$
\[ T_{n,h}\le \Delta_1(\Sigma_X,h)+{\mathcal O}_P(\Delta_2(\Sigma_X,h)(nh)^{-1/2})=\Delta_1(\Sigma_X,h)+{\mathcal O}_P(n^{-1/2}),\]
where we use $\Delta_2(\Sigma_X,h)\le h^{1/2}\abs{\Sigma_X}_{B^{1/2}_{2,\infty}}$ and $\Sigma_X\in B^{1/2}_{2,\infty}$ almost surely. Therefore Proposition \ref{PropNQV} and $h+(nh)^{-1}=o(h')$ yield on ${\mathcal H}_0$
\[ T_{n,h}\le \tfrac{8}{15}h^{1/2}NV^{(1)}+{\mathcal O}_P(h+n^{-1/2})= \tfrac{8}{15}h^{1/2}\big(NV^{(1)}+{ o}_P((h')^{1/2})\big).\]
By the Central Limit Theorem \ref{ThmNQVCLT}, the consistency results of Proposition \ref{PropVarhat} and Slutsky's lemma we conclude
\[ \Big(3h' \big(\widehat{NV}_{h',n}^{(2)}-\widehat{BNV}^{(1)}_{h',n}\big)\Big)^{-1/2}\Big(\widehat{NV}^{(1)}_{h',n}-NV^1\Big) \xrightarrow{d} N(0,1).
\]
Another application of Slutsky's lemma implies that
\[ \liminf_{n\to\infty}\PP\Big(\Big\{\Big(3h' \big(\widehat{NV}_{h',n}^{(2)}-\widehat{BNV}^{(1)}_{h',n}\big)\Big)^{-1/2}\Big(\tfrac{15}{8h^{1/2}}T_{n,h}-\widehat{NV}^{(1)}_{h',n}\Big)\le q_{1-\alpha;N(0,1)}\Big\}\cap{\mathcal H}_0\Big)\]
is at least $1-\alpha$. This gives the result for the test $\phi_{\alpha,n,h,h'}$.
\end{proof}

In \citet[Theorem 2]{FanWang2008} it is shown that a standard kernel estimator $\hat\Sigma_X(t)$ of the spot covariance $\Sigma_X(t)$ is consistent uniformly over $t\in[0,1]$. Moreover, the convergence rate $n^{-1/4}$ up to log factors is derived together with a limiting distribution. Since the eigenvalue map $\Sigma\mapsto \lambda_r(\Sigma)$ is Lipschitz continuous (with respect to the spectral norm), the same uniform rate of convergence holds for $\lambda_r(\hat\Sigma_X(t))$ such that
\begin{equation}\label{Eqlambdahat}
\hat{\underline\lambda}_r:=\inf_{t\in[0,1]}\lambda_r(\hat\Sigma_X(t))=\inf_{t\in[0,1]}\lambda_r(\Sigma_X(t))+{\mathcal O}_P(n^{-1/6}).
\end{equation}
These properties can also be derived for our blockwise realised covariance matrix $\hat\Sigma^{kh}$ as an estimator of $\Sigma_X(t)$ for $t\in I_k$ when $h\thicksim n^{-1/2}$.
This enables completely data-driven critical values also in the spectral gap case. The next result is proved in Appendix \ref{AppStochVol}.

\begin{corollary}\label{Corkappahat2}
Assume Model \eqref{EqStochVol} under Assumption \ref{AssStochVol} for $p=2$. Consider block sizes $h,h'\to 0$ with $h=o(h')$, $h'n^{1/3}\to\infty$ and $nhh'\to\infty$. Let $\hat{\underline\lambda}_r$ be a spectral gap estimator satisfying \eqref{Eqlambdahat}.
Then the test $\phi_\alpha=\phi_{\alpha,n,h,h'}$ from \eqref{Eqphi0} with fixed $\alpha\in(0,1)$ and critical value
\[ \kappa_\alpha=\kappa_{\alpha,n,h,h'}=\frac {h}{3\hat{\underline\lambda}_r}\Big(\widehat{NV}_{h',n}^{(2)}+ \Big(3h'\big(\widehat{NV}_{h',n}^{(4)}-\widehat{BNV}_{h',n}^{(2)}\big)\Big)^{1/2}q_{1-\alpha;N(0,1)}\Big)\]
has asymptotic level $\alpha$ on the (possibly random) null hypothesis ${\mathcal H}_0^{gap}=\{\sup_{t\in[0,1]}\rank(\Sigma_X(t))\le r,\,\inf_{t\in[0,1]}\lambda_r(\Sigma_X(t))>0\}$:
\[ \limsup_{n\to\infty}\PP\Big(\{\phi_{\alpha,n,h,h'}=1\}\cap{\mathcal H}_0^{gap}\Big)\le\alpha.\]
\end{corollary}

\begin{figure}[t]
\centering
\includegraphics[width=0.32\textwidth]{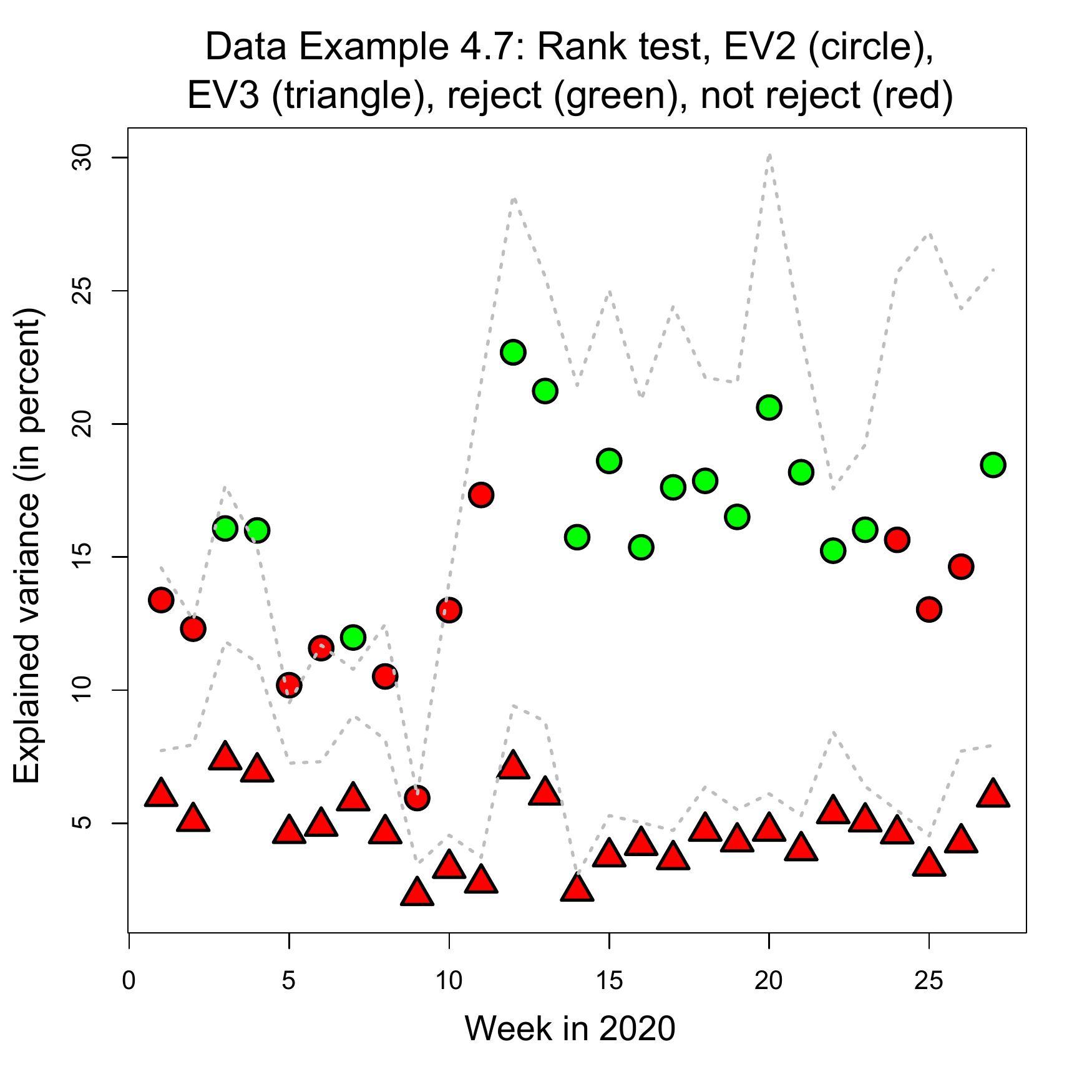}
 \includegraphics[width=0.32\textwidth]{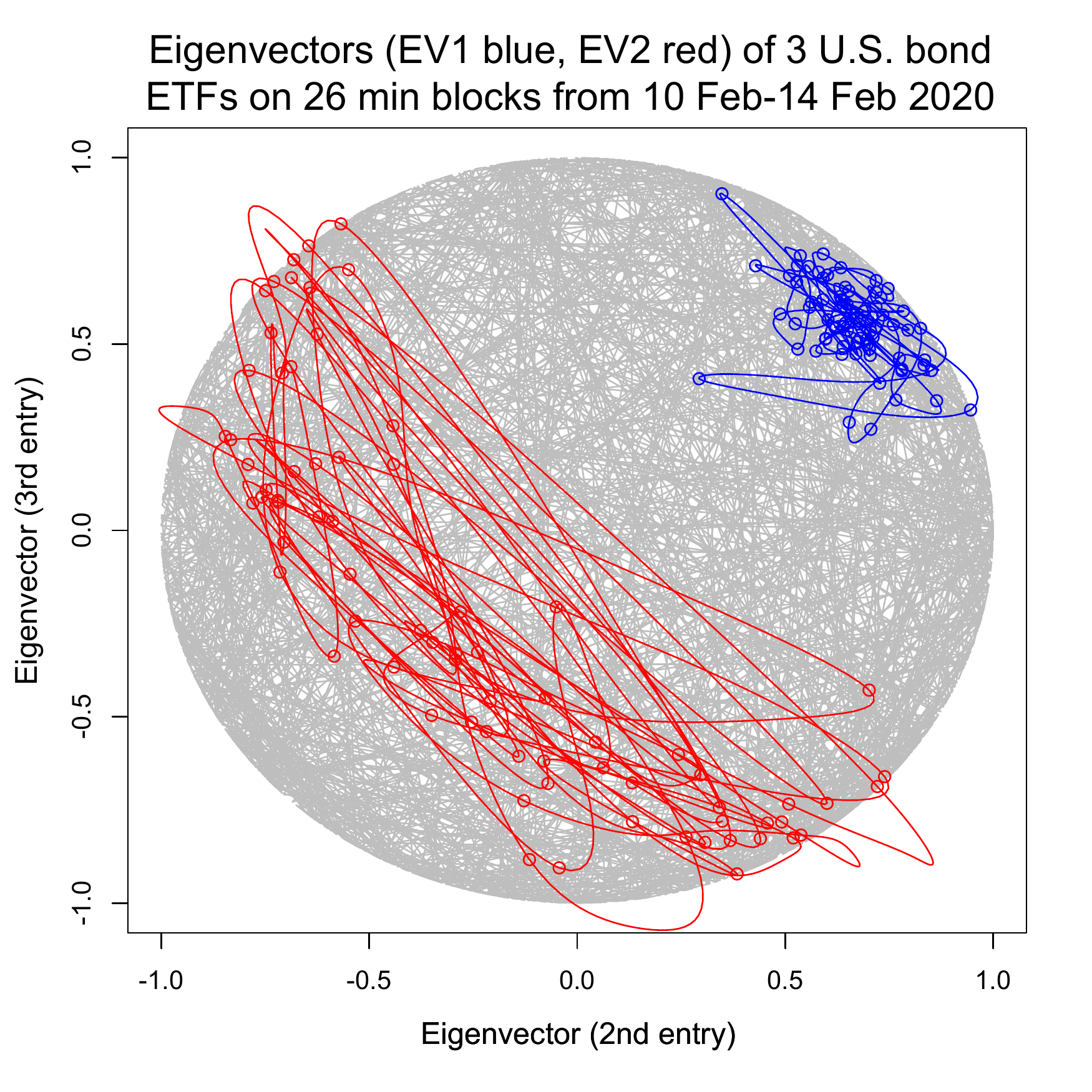}
\includegraphics[width=0.32\textwidth]{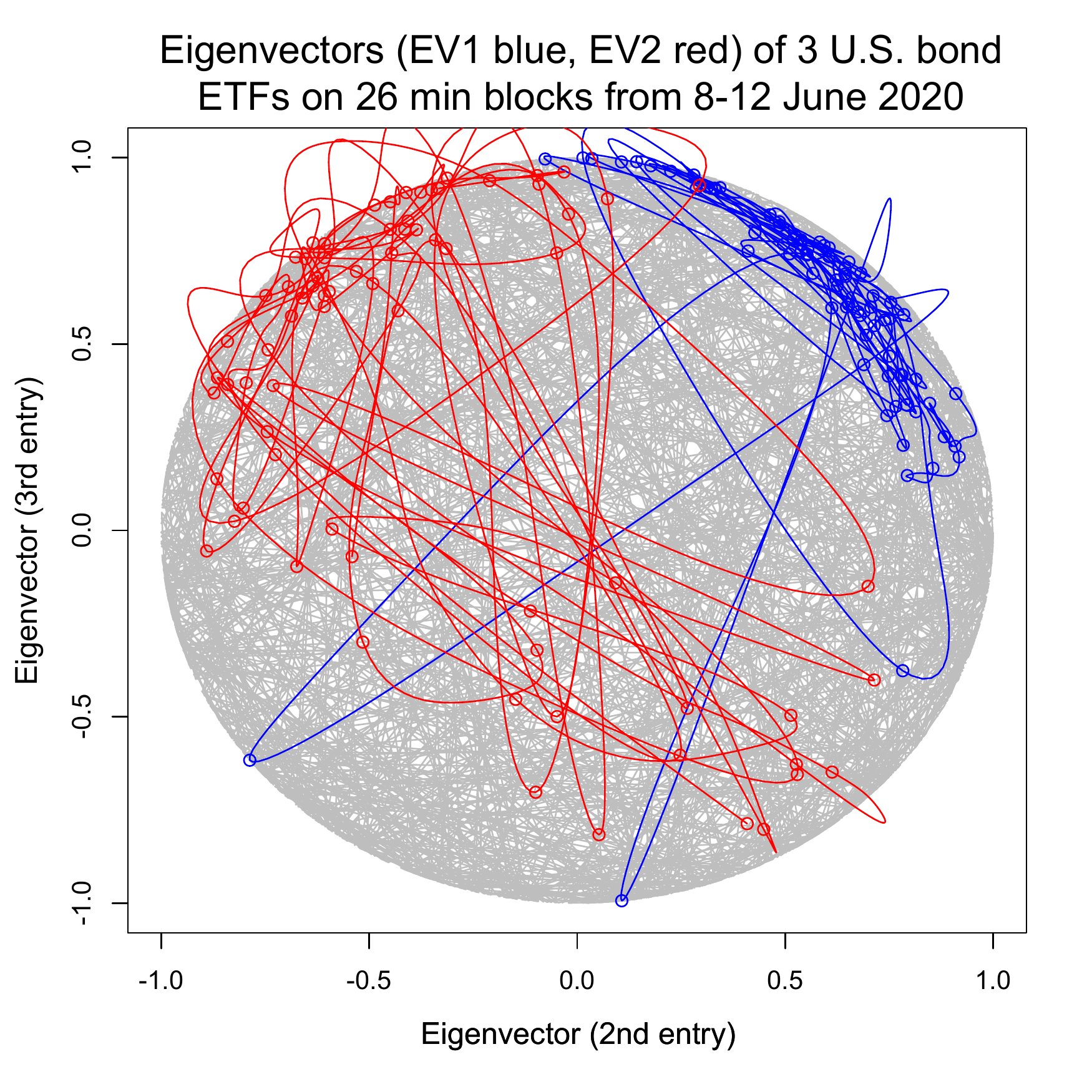}
 \caption{U.S.\ government bond Examples \ref{ExEmpirical} and \ref{ExEmpirical2}. Left: Test results for each week. Center: Movements of eigenvectors in week 7. Right: Movements of eigenvectors in week 24.}\label{FigBondData}
\end{figure}

\begin{example}\label{ExEmpirical2}
We apply the test with estimated critical values to the bond data already studied in Example \ref{ExEmpirical} above.
Choose the block length $h=0.013$ (26 minutes), that is 15 blocks per day and $h^{-1}=75$ blocks per week. On those 75 blocks we calculate the test statistic \eqref{Eqphi0} for $r=1$ and $r=2$. Critical values are taken from Corollary \ref{Corkappahat2} with estimators $\widehat{NV}_{h',n}$ and $\widehat{BNV}_{h',n}$ taken on blocks of size $h'=0.033$ (65 minutes). We use $\hat{\underline\lambda}_r=\min \lambda_{r+1}(\hat\Sigma^{kh})$ and choose $\alpha=5\%$.

Figure \ref{FigBondData}(left) displays the test results along with the explained variance  for each of the first 27 weeks in 2020.  Circles and triangles show the explained variances \eqref{EqExplVar} for the second and third eigenvalues, respectively. The colours indicate whether the test accepts (red) or rejects (green) the hypothesis of rank 1 (circles) or rank 2 (triangles). We see that in some weeks rank 1 is accepted, while always a maximal rank 2 is accepted versus rank 3 alternatives. The dashed lines visualise the corresponding explained variances \eqref{EqExplVar} when weekly integrated covariance matrices $\Sigma_X$ ($h=1$) are used. Notice that the time period includes the Corona pandemic and major disruptions of financial markets from March (week 10) onwards.

We highlight two specific weeks. We find that the test rejects rank 1 in week 7 (Feb.\ 10-14) but fails to reject in week 24 (June 8-12). Both weeks seem to exhibit with about 12\% and 15\% similarly large explained variance of a second principal component. Figure \ref{FigBondData}(center and right) reveal the source of the different test decisions by plotting the normalized (length one, first entry positive) first (blue) and second (red) eigenvectors of $\hat\Sigma_X^{kh}$ on the upper hemisphere of the unit ball. The movement of eigenvectors is much stronger in week 24 than in week 7. The larger the movements the larger are the calibrated critical value. In other words, faster dynamics may lead to a larger test statistics $T_{n,h}$ under the null, which is why rank 1 is not rejected in week 24. The example shows how periods of high (co)volatility coincide with periods of large movements in eigenvectors and how this can affect the inference on the rank.

Our data analysis does not give empirical evidence that the U.S.\ term structure for maturities between 3 and 20 years exceeds rank 2 in the first six months of 2020.
\end{example}

\begin{appendix}

\section{Further results and proofs}

\subsection{Matrix deviation results}\label{secMatrixConc}

We provide explicit deviation bounds for eigenvalues under averaging.

\begin{proof}[Proof of Proposition \ref{PropPerturb}]
Set $\bar S=\frac1{\abs{I}}\int_I S(t)dt$ and denote by $\bar\lambda_1\ge\cdots\ge\bar\lambda_d$  the ordered eigenvalues (with multiplicities) of $\bar S$. The standard eigenvalue-norm bound implies for any $t\in I$
\[ \bar\lambda_{r+1}\le \lambda_{r+1}(S(t))+\norm{\bar S-S(t)}=\norm{\bar S-S(t)}
\]
because of $\rank(S(t))\le r$. In particular, it is smaller than the mean of $\norm{\bar S-S(t)}$ over $t$ and we conclude by the convexity of the norm
\[  \bar\lambda_{r+1}\le \frac1{\abs{I}} \int_I \norm{\bar S-S(t)}\,dt\le \frac1{\abs{I}^2}\int_{I\times I} \norm{S(s)-S(t)}\,dt\,ds=\Delta_1(S,I).
\]
Since for $\underline\lambda_r<2 \Delta_1(S,I)$ the minimum in the asserted inequality is attained at $\Delta_1(S,I)$ in view of $\Delta_2(S,I)\ge \Delta_1(S,I)$, it remains to show $\bar\lambda_{r+1} \le \tfrac{2}{\underline\lambda_r}\Delta_2(S,I)^2$ under the assumption $\underline\lambda_r\ge 2\Delta_1(S,I)$.

We use the spectral decomposition  $\bar S =\sum_{j=1}^d\bar\lambda_j\bar P_j$ with the rank-one projections $\bar P_j$ onto the eigenspaces of $\bar S$ corresponding to $\bar\lambda_j$. We apply the resolvent identity
\[\bar P_j=\bar P_j(\bar S- \bar\lambda_j +\bar\lambda_j-S(t))(\bar\lambda_j-S(t))^{-1}=\bar P_j(\bar S-S(t))(\bar\lambda_j-S(t))^{-1},\]
assuming that $\bar\lambda_{j}$ is not an eigenvalue of $S(t)$, and obtain
\begin{align*}
\bar\lambda_{j}&=\scapro{\bar S-S(t)}{\bar P_{j}}_{HS}+\scapro{S(t)}{\bar P_{j}(\bar S-S(t))(\bar\lambda_{j}-S(t))^{-1}}_{HS}.
\end{align*}
Taking transposes, we also have $\bar P_j=(\bar\lambda_j-S(t))^{-1}(\bar S-S(t))\bar P_j$ and we can further expand
\begin{align}\label{EqQuadrPert}
\bar\lambda_{j}&=\scapro{\bar S-S(t)}{\bar P_{j}}_{HS}+\scapro{S(t)}{(\bar\lambda_{j}-S(t))^{-1}(\bar S-S(t))\bar P_{j}(\bar S-S(t))(\bar\lambda_{j}-S(t))^{-1}}_{HS}.
\end{align}
Therefore by functional calculus with $f(\lambda):=\lambda(\lambda-\bar\lambda_{r+1})^{-2}$ we have for $j=r+1$
\[ \bar\lambda_{r+1}=\scapro{\bar S-S(t)}{\bar P_{r+1}}_{HS}+\scapro{f(S(t))}{(\bar S-S(t))\bar P_{r+1}(\bar S-S(t))}_{HS}.
\]
We know $\bar \lambda_{r+1}\le \Delta_1(S,I)\le  \frac12\underline\lambda_r$  by the first part and the assumption on $\underline\lambda_r$. Hence, we conclude $\bar \lambda_{r+1}<\underline\lambda_r$ and $\bar \lambda_{r+1}$ is indeed not an eigenvalue of $S(t)$ for any $t\in I$.
By integrating over $I$, the linear term vanishes and
\begin{align*}
\bar\lambda_{r+1} &=\frac1{\abs{I}}\int_I \trace\Big(f(S(t))(\bar S-S(t))\bar P_{r+1}(\bar S-S(t))\Big)\,dt\\
 &=\frac1{\abs{I}} \trace\Big(\int_I(\bar S-S(t))f(S(t))(\bar S-S(t))\,dt\bar P_{r+1}\Big).
\end{align*}
Since $\rank(S(t))\le r$ and $\lambda_r(S(t))\ge\underline\lambda_r$, we have $\norm{f(S(t))}\le \underline\lambda_r(\underline\lambda_r-\bar\lambda_{r+1})^{-2}$ and we can bound
\[ \bar\lambda_{r+1}\le \frac{\underline\lambda_r}{(\underline\lambda_r-\bar \lambda_{r+1})^2}\bnorm{\frac1{\abs{I}}\int_I (\bar S-S(t))^2dt} \trace(\bar P_{r+1}).
\]
Finally, use
\[ \bnorm{\frac1{\abs{I}}\int_I (\bar S-S(t))^2dt}=\bnorm{\frac1{\abs{I}}\int_I S(t)^2dt-\Big(\frac1{\abs{I}}\int_I S(t)\,dt\Big)^2}=\frac12\Delta_2(S,I)^2
\]
as well as $\underline\lambda_r-\bar \lambda_{r+1}\ge \frac12\underline\lambda_r$ and $\trace(\bar P_{r+1})=1$ to obtain the  assertion.
\end{proof}

\begin{proposition}\label{PropNormPerturb2}
In the setting of Proposition \ref{PropPerturb} denote by $S(t)_{>r}\in\R^{(d-r)\times(d-r)}$ the  minor of $S(t)$ on the space $V_{>r}$ generated by the $d-r$ smallest eigenvalues of $\bar S$, comparable to \eqref{EqMinor} below.
\begin{enumerate}
\item We always have
\[ \norm{S(t)_{>r}}\le \lambda_{r+1}(\bar S)+\norm{S(t)-\bar S}.\]
\item If $\lambda_r(S(t))> \lambda_{r+1}(\bar S)$ holds, then also
\[ \norm{S(t)_{>r}}\le \frac{r\lambda_r(S(t))}{(\lambda_r(S(t))-\lambda_{r+1}(\bar S))^2}\norm{S(t)-\bar S}^2.
\]
\item For $I'\subset I$, $p\ge 1$ and $\underline\lambda_r=\inf_{t\in I'}\lambda_r(S(t))$ we always have in terms of $p$-variations
\[ \bnorm{\frac1{\abs{I'}}\int_{I'}S(t)_{>r}dt}\le \Big(\big(1+(\abs{I}/\abs{I'})^{1/p}\big)\Delta_p(S,I)\Big)\wedge \Big((r+4)(\abs{I}/\abs{I'})^{1/p}\frac{\Delta_{2p}(S,I)^2}{\underline\lambda_r}\Big).
\]
\end{enumerate}
\end{proposition}

\begin{proof}
Write $\bar\lambda_j=\lambda_j(\bar S)$.
For (a) we argue by  triangle inequality and $\norm{\bar P_{>r}}\le 1$:
\[ \norm{S(t)_{>r}}=\norm{\bar P_{>r}S(t)\bar P_{>r}}\le   \norm{\bar P_{>r}\bar S\bar P_{>r}}+ \norm{\bar P_{>r}(S(t)-\bar S)\bar P_{>r}}
\le \bar\lambda_{r+1}+\norm{S(t)-\bar S}.
\]

For (b) we proceed as for the derivation of \eqref{EqQuadrPert}. For $j>r$ consider $f_j(\lambda):=\lambda(\lambda-\bar\lambda_j)^{-2}$ and use  $f_j(S(t))=P_{\ker(S(t))^\perp}f_j(S(t))P_{\ker(S(t))^\perp}$ with the orthogonal projection $P_{\ker(S(t))^\perp}$ onto the orthogonal complement of the kernel of $S(t)$ and equivalently onto the range of $S(t)$. Then
\begin{align*}
\scapro{S(t)}{\bar P_j}_{HS} &= \scapro{S(t)}{(\bar\lambda_{j}-S(t))^{-1}(\bar S-S(t))\bar P_{j}(\bar S-S(t))(\bar\lambda_{j}-S(t))^{-1}}_{HS}\\
 &= \trace\Big(f_j(S(t))P_{\ker(S(t))^\perp}(\bar S-S(t))\bar P_{j}(\bar S-S(t))P_{\ker(S(t))^\perp}\Big)\\
&\le \norm{f_j(S(t))}\trace\Big(P_{\ker(S(t))^\perp}(\bar S-S(t))\bar P_j(\bar S-S(t))P_{\ker(S(t))^\perp}\Big).
\end{align*}
Because of $\lambda_r(S(t))>\bar\lambda_{r+1}\ge\bar\lambda_j$ and $\rank(S(t))=r$ we can bound $\norm{f_j(S(t))}\le \lambda_r(S(t))(\lambda_r(S(t))-\bar\lambda_{r+1})^{-2}$. We expand $\bar P_{>r}=\sum_{j>r}\bar P_j$ and arrive at
\begin{align*}
\norm{S(t)_{>r}} &\le \trace(\bar P_{>r}S(t)\bar P_{>r})= \sum_{j=r+1}^d \scapro{S(t)}{\bar P_j}_{HS}\\
&\le \frac{\lambda_r(S(t))}{(\lambda_r(S(t))-\bar\lambda_{r+1})^2}\trace\Big(P_{\ker(S(t))^\perp}(\bar S-S(t))\bar P_{>r}(\bar S-S(t))P_{\ker(S(t))^\perp}\Big)\\
&\le \frac{\lambda_r(S(t))}{(\lambda_r(S(t))-\bar\lambda_{r+1})^2}\trace\big(P_{\ker(S(t))^\perp}\big)\norm{\bar S-S(t)}^2\norm{\bar P_{>r}}\\
&= \frac{\lambda_r(S(t))}{(\lambda_r(S(t))-\bar\lambda_{r+1})^2}r\norm{\bar S-S(t)}^2.
\end{align*}

To obtain (c), we integrate the previous bounds. From (a) and Proposition \ref{PropPerturb} we obtain  by the convexity of the norm and H\"older's inequality
\begin{align*}
\int_{I'}\norm{S(t)_{>r}}\,dt &\le \abs{I'}\Delta_1(S,I)+\int_{I'}\frac1{\abs{I}}\int_I \norm{S(t)-S(s)}\,dtds\\
&\le \abs{I'}\Delta_1(S,I)+\abs{I'}^{1-1/p}\Big(\int_{I'}\frac1{\abs{I}}\int_I \norm{S(t)-\bar S}^pdt\Big)^{1/p}\\
&\le \abs{I'}\Delta_1(S,I)+\abs{I'}^{1-1/p}\abs{I}^{1/p}\Delta_p(S,I),
\end{align*}
where we bounded the integral over $I'$ by the integral over $I$ in the last step. Bounding $\Delta_1(S,I)\le \Delta_p(S,I)$ yields
\[ \bnorm{\frac1{\abs{I'}}\int_{I'}S(t)_{>r}\,dt} \le \big(1+(\abs{I}/\abs{I'})^{1/p}\big)\Delta_p(S,I).\]
This gives the assertion for $\underline\lambda_r\le \frac{r+4}{2}\Delta_1(S,I)\le \frac{r+4}{2}\Delta_{2p}(S,I)$ because the minimum in the asserted bound is then attained at the first argument. Otherwise, $\bar\lambda_{r+1}\le \Delta_1(S,I)\le \frac2{r+4}\underline\lambda_r$ holds by Proposition \ref{PropPerturb}. Then the same arguments as before, but using part (b) give
\begin{align*}
\bnorm{\frac1{\abs{I'}}\int_{I'}S(t)_{>r}\,dt }
&\le (\abs{I}/\abs{I'})^{1/p}\frac{r\underline\lambda_r }{(\underline\lambda_r-\bar\lambda_{r+1})^2}\Delta_{2p}(S,I)^2\\
&\le \frac{r(r+4)^2}{(r+2)^2}(\abs{I}/\abs{I'})^{1/p}\frac{\Delta_{2p}(S,I)^2}{\underline\lambda_r}.
\end{align*}
The claim therefore follows with $r(r+4)\le (r+2)^2$.
\end{proof}

\subsection{Upper matrix concentration bounds}

Concentration results for the maximal eigenvalue of sums over independent, but non-identically distributed Wishart matrices are established. We follow the proof of the matrix Bernstein inequalities in \citet{tropp2012}, but argue differently from the subexponential case there because we face non-commuting matrices.

\begin{theorem}\label{ThmBernstein}
Let $Y_j\sim N(0,A_j)$, $j=1,\ldots,J$, be independent Gaussian random vectors for some $A_j\in\R_{spd}^{d\times d}$. Let $A:=\sum_{j= 1}^J A_j$ and
\begin{align*}
\sigma^2 &:=\lambda_{max}\Big(\sum_{j= 1}^J(\trace(A_j) A_j+2A_j^2)\Big),\quad
 R:=\max_{j=1,\ldots,J}\big(\trace(A_j)+4\norm{A_j}\big).
\end{align*}
Then we have  the upper tail bound
\begin{align*}
&\PP\Big(\lambda_{max}\Big(\sum_{j= 1}^J Y_jY_j^\top-A\Big)\ge t\Big)\le d\exp\Big(-\frac{t^2}{2\sigma^2+2Rt}\Big)
\end{align*}
and the expectation bound
\begin{align*}
&\E\Big[\lambda_{max}\Big(\sum_{j= 1}^J Y_jY_j^\top-A\Big)\Big]\le \sigma\big(2\sqrt{\log d}+1\big)+4R\big(\log(d)+1\big).
\end{align*}
\end{theorem}

\begin{proof}
Let us write $S_j=Y_jY_j^\top-A_j$ and note $\E[S_j]=0$, $S_j\ge -A_j$. We start with the master tail bound from the trace exponential \citep[Theorem 3.6]{tropp2012} and obtain
\begin{align*}
\PP\Big(\lambda_{max}\Big(\sum_{j\ge 1} Y_jY_j^\top-A\Big)\ge t\Big) &= \PP\Big(\lambda_{max}\Big(\sum_{j\ge 1} S_j\Big)\ge t\Big)\\
&\le \inf_{\theta>0} \trace\Big(\exp\Big(-\theta t I_d+\sum_{j\ge 1} \log\big(\E\big[e^{\theta S_j}\big]\big)\Big)\Big).
\end{align*}
The function $f(x)=e^x-1-x-x^2/2$ is increasing on $\R$ and $e^x-1-x\le \frac12 x^2e^x$ holds for $x\ge 0$, which is easily checked by a power series expansion. Functional calculus and $S_j\le Y_jY_j^\top$ then show
 \[ e^{\theta S_j}-I_d-\theta S_j\le \tfrac12(\theta S_j)^2+f(Y_jY_j^\top)\le \frac{\theta^2}2\Big(S_j^2-(Y_jY_j^\top)^2+(Y_jY_j^\top)^2e^{\theta Y_jY_j^\top}\Big).
\]
We use $\log(I_d+A)\le A$ and $\E[S_j]=0$ to arrive at
\begin{align*}
\log\big(\E[e^{\theta S_j}]\big)&\le \E[e^{\theta S_j}-I_d-\theta S_j]\\
&\le \frac{\theta^2}{2}\Big(-A_j^2+\E\big[(Y_jY_j^\top)^2\exp(\theta Y_jY_j^\top)\big]\Big)\\
&\le \frac{\theta^2}{2}\E\Big[(Y_jY_j^\top)^2\exp(\theta Y_jY_j^\top)\Big].
\end{align*}
Observe that $Y_jY_j^\top$ is a rank-one matrix with non-zero eigenvalue $\norm{Y_j}^2$ and corresponding eigenprojector $Y_jY_j^\top/\norm{Y_j}^2$. By functional calculus we thus have
\begin{align*}
(Y_jY_j^\top)^2\exp(\theta Y_jY_j^\top) &= \norm{Y_j}^2\exp(\theta \norm{Y_j}^2)Y_jY_j^\top.
\end{align*}
Writing $Y_j=(Y_{j,1},\ldots,Y_{j,d})^\top$ in the basis of eigenvectors of $A_j$, the coordinates $Y_{j,i}$ are independent Gaussian $N(0,\lambda_i(A_j))$-distributed, in particular $Y_{j,i}Y_{j,\ell}\stackrel{d}{=}-Y_{j,i}Y_{j,\ell}$ holds for $i\not=\ell$. This implies
\[  \E\big[(Y_jY_j^\top)^2\exp(\theta Y_jY_j^\top)\big]=\diag\Big(\E\big[\norm{Y_j}^2\exp(\theta \norm{Y_j}^2)Y_{j,i}^2\big]\Big)_{i=1,\ldots,d}.
\]
For $Z\sim N(0,\lambda)$ simple $\Gamma$-density identities yield for $\theta\le (2\lambda)^{-1}$
\[ \E\big[Z^{2m}e^{\theta Z^2}]=\frac{\Gamma(m+1/2)(2\lambda)^m}{\Gamma(1/2)(1-2\theta\lambda)^{m+1/2}}.\]
Using independence of the coordinates and this result for $m=0,1,2$, we find for $\theta\le (2\norm{A_j})^{-1}$
\begin{align*}
&\diag\Big(\E\big[\norm{Y_j}^2\exp(\theta \norm{Y_j}^2)Y_{j,i}^2\big]\Big)_{i=1,\ldots,d} =\Big(\prod_{i=1}^d(1-2\theta\lambda_i(A_j))^{-1/2}\Big) \\ &\quad\times\diag\Big(\frac{2\lambda_i^2(A_j)}{(1-2\theta\lambda_i(A_j))^2}
+\frac{\lambda_i(A_j)}{1-2\theta\lambda_i(A_j)}\sum_{\ell=1}^d\frac{\lambda_\ell(A_j)}{1-2\theta\lambda_\ell(A_j)}\Big)_{i=1,\ldots,d}\\
&=\det(I-2\theta A_j)^{-1/2}\Big(2A_j^2(I-2\theta A_j)^{-2}
+A_j(I-2\theta A_j)^{-1}\trace(A_j(I-2\theta A_j)^{-1})\Big).
\end{align*}
We derive a simpler bound, using $\prod_i(1-2a_i)^{-1/2}\le (1-\sum_ia_i)^{-1}$ for $0\le a_i\le 1/2$ and $\trace(A_j(1-2\theta A_j)^{-1})\le \trace(A_j)(1-2\theta\norm{A_j})^{-1}$:
\[ \E\big[(Y_jY_j^\top)^2\exp(\theta Y_jY_j^\top)\big]\le (1-\theta(\trace(A_j)+4\norm{A_j}))^{-1}\big(2A_j^2+A_j\trace(A_j)\big).\]
So far, we have thus established for $\theta\in(0,R^{-1})$
\[ \sum_{j=1}^J \log\big(\E\big[e^{\theta S_j}\big]\big)\le \frac{\theta}{2(\theta^{-1}-R)} \sum_{j=1}^J\big(\trace(A_j) A_j+2A_j^2\big).
\]
Inserting $\theta=t(\sigma^2+tR)^{-1}$ and applying $\trace(\exp(M))\le de^{\lambda_{max}(M)}$ for symmetric matrices $M$, the upper tail bound follows.


The expectation bound follows from integrating
\[\PP\Big(\lambda_{max}\Big(\sum_{j\ge 1} Y_jY_j^\top-A\Big)\ge t\Big)\le 1\wedge\Big(d\exp\Big(-\frac{t^2}{4\sigma^2}\Big)+ d\exp\Big(-\frac{t}{4R}\Big)\Big)
\]
over $t\in[0,\infty)$.
\end{proof}

\begin{remark}
It is easy to check that $\sigma^2=\lambda_{max}(\sum_j\E[(Y_jY_j^\top)^2])$ holds as usual for matrix Bernstein inequalities.
 In dimension $d=1$ the so far best known deviation bound by \citet{laurent2000} is
\[ \PP\Big(\sum_{j= 1}^J (Y_j^2-A_j)>t\Big)\le \exp\Big(-\frac{t^2}{4\sum_{j}A_j^2+4t\max_{j}A_j}\Big).\]
We observe that we merely lose by a linear dependence of the constants in the dimension for the matrix concentration bound. We could argue as \citet{laurent2000} and win in the constants for the isotropic case $A_j=I_d$, but in general we face the problem that $Y_jY_j^\top$ and $A_j$ do not commute.
\end{remark}

We combine the matrix Bernstein inequality with a Gaussian concentration argument to obtain a subexponential deviation inequality for a triangular scheme of maximal eigenvalues.

\begin{theorem}\label{ThmBernsteinTriangular}
Let $Y_{jk}\sim N(0,A_{jk})$, $j=1,\ldots,J$, $k=1,\ldots,K$, be independent Gaussian random vectors for some  $A_{jk}\in\R_{spd}^{d\times d}$. Let
\begin{align*}
\sigma_k^2 &:=\lambda_{max}\Big(\sum_{j= 1}^J(\trace(A_{jk}) A_{jk}+2A_{jk}^2)\Big),\quad
 R_k:=\max_{j=1,\ldots,J} (\trace(A_{jk})+4\norm{A_{jk}}).
\end{align*}
and $A_k:=\sum_{j=1}^JA_{jk}$. Then we have with probability at least $1-e^{-t}$ for any $\delta>0$
\begin{align*}
\sum_{k=1}^K\lambda_{max}\Big(\sum_{j=1}^J Y_{jk}Y_{jk}^\top\Big)&\le  (1+\delta)\sum_{k=1}^K\Big(\lambda_{max}(A_k)+\sigma_k(2\sqrt{\log d}+1)+4R_k(\log(d)+1)\Big)\\
&\quad + 2(1+\delta^{-1})\max_{j,k}\norm{A_{jk}} t.
\end{align*}
\end{theorem}

\begin{proof}
The expectation bound of Theorem \ref{ThmBernstein} together with the inequality $\lambda_{max}(B)\le \lambda_{max}(B-A)+\lambda_{max}(A)$ yields
\[ \E\Big[\sum_{k=1}^K\lambda_{max}\Big(\sum_{j=1}^J Y_{jk}Y_{jk}^\top\Big)\Big]\le \sum_{k=1}^K\Big(\lambda_{max}(A_k)+\sigma_k(2\sqrt{\log d}+1)+4R_k(\log(d)+1)\Big).
\]
We write $Y_{jk}=A_{jk}^{1/2}Z_{jk}$ with independent $Z_{jk}\sim N(0,I_d)$  and use the Gaussian concentration \citep[Equation (1.6)]{LedouxTal1991} to deduce for $\kappa>0$ that
\[ \PP\Big(\Big(\sum_{k=1}^K\lambda_{max}\Big(\sum_{j= 1}^J A_{jk}^{1/2}Z_{jk}Z_{jk}^\top A_{jk}^{1/2}\Big)\Big)^{1/2}\ge \E\Big[\Big(\sum_{k=1}^K\lambda_{max}\Big(\sum_{j= 1}^J Y_{jk}Y_{jk}^\top\Big)\Big)^{1/2}\Big]+\kappa L\Big)
\]
is at most $e^{-\kappa^2/2}$, where $L$ is the Lipschitz constant of the functional $F((z_{jk})_{j,k}):=(\sum_{k=1}^K\lambda_{max}(\sum_{j= 1}^J A_{jk}^{1/2}z_{jk}z_{jk}^\top A_{jk}^{1/2}))^{1/2}$. $F$ defines a norm on $(\R^d)^{J\times K}$ so that
\begin{align*}
L^2&=\sup_{(z_{jk})_{j,k}\not=0} \frac{\sum_{k=1}^K\sup_{\norm{w_k}\le 1}\sum_{j= 1}^J \scapro{A_{jk}^{1/2}z_{jk}}{w_k}^2} {\sum_{k=1}^K\sum_{j=1}^J \norm{z_{jk}}^2}\\
&\le \sup_{(z_{jk})_{j,k}\not=0} \frac{\sum_{k=1}^K\sum_{j= 1}^J \norm{{A_{jk}^{1/2}z_{jk}}}^2} {\sum_{k=1}^K\sum_{j=1}^J \norm{z_{jk}}^2}= \max_{j,k}\norm{A_{jk}}.
\end{align*}
Taking squares and applying the Cauchy-Schwarz inequality for the expectation as well as $(A+B)^2\le (1+\delta)A^2+(1+\delta^{-1})B^2$ for $A,B,\delta>0$, we thus find with probability at least $1-e^{-\kappa^2/2}$
\begin{align*}
\sum_{k=1}^K\lambda_{max}\Big(\sum_{j= 1}^J Y_{jk}Y_{jk}^\top\Big)&\le  (1+\delta)\E\Big[\sum_{k=1}^K\lambda_{max}\Big(\sum_{j=1}^J Y_{jk}Y_{jk}^\top\Big)\Big]+ (1+\delta^{-1})\kappa^2\max_{j,k}\norm{A_{jk}}.
\end{align*}
It remains to insert the expectation bound and to set $t=\kappa^2/2$.
\end{proof}

\subsection{Lower matrix concentration bounds}

Lower tail bounds for the minimal eigenvalue are required to establish the asymptotic power of the tests under the alternative and need not be precise in the constants for our needs. They must, however, be tight in the sense that they prescribe correctly the probability for minimal eigenvalues approaching zero. The standard Chernoff bounds for bounded random matrices \citep{tropp2012} or the concentration bounds for empirical covariance matrices \citep{koltchinskii2015} still produce positive tail bounds at zero, which only for increasing sample size become negligible. We also need to cope with constant sample sizes on each block while averaging over an increasing number of blocks.

The proof is surprisingly intricate because the covariance matrices  need not be jointly diagonalisable. First, we establish a stochastic dominance property for sums of Wishart matrices via explicit density calculations when each population covariance matrix is larger than a fixed covariance matrix. In a second step we use the entropy of the Grassmanian manifold to exhibit such a fixed covariance matrix for a fraction of all summands. We are not aware of any more direct argument. In particular, the number of summands required is pretty large and a tighter bound would certainly be desirable.
Let us start with determining explicitly the eigenvalue density for sums of different Wishart matrices.

\begin{lemma}\label{LemWishartGen}
Let $d,J\in\N$ with $J\ge d$ and consider $Y\in\R^{d\times J}$, a centred Gaussian matrix with $\vek(Y)\sim N(0,A)$ and an invertible covariance matrix $A\in\R^{dJ\times dJ}$. Then the eigenvalues $\lambda_1\ge\cdots\ge \lambda_d\ge 0$ of $YY^\top$ have  the joint Lebesgue density
\begin{align*}
f^\Lambda(\lambda_1,\ldots,\lambda_{d})&=c\det(A)^{-1/2}\prod_{j=1}^{d}\lambda_j^{(J-d-1)/2}\prod_{i<j}(\lambda_i-\lambda_j)\times\\
&
\iint \exp\Big(-\frac12\scapro{A^{-1}\vek(O' L O)}{\vek(O'LO)}_{\R^{dJ}}\Big)\,\HH_J(dO)\HH_d(dO')
\end{align*}
with $c=\pi^{d^2/2}2^{-dJ/2}\Gamma(d/2)^{-1}\Gamma(J/2)^{-1}$ and $L=(\diag(\lambda_1^{1/2},\ldots,\lambda_{d}^{1/2}),0_{d\times (J-d)})\in\R^{d\times J}$ in terms of the zero matrix $0_{d\times (J-d)}\in\R^{d\times(J-d)}$. We write $\int g(O)\HH_m(dO)$ when integrating $g:O(m)\to\R$ with respect to the normalised Haar measure $\HH_m$ on the orthogonal matrix group $O(m)$.
\end{lemma}

\begin{proof}
Extending the density result by \citet{James1960} to the non-i.i.d. setting, we consider the random matrix $O'YO$ with orthogonal transformations  $(O',O)\sim \HH_d\otimes \HH_J$ independent of $Y$. The idea is that $O'YO(O'YO)^\top$ and $YY^\top$ share the same eigenvalues. The Lebesgue density of $O'YO$ is given by the Gaussian mixture
\begin{align*}
f^{O'YO}(y)&=(2\pi)^{-dJ/2} \det(A)^{-1/2}\times\\
&\quad  \iint \exp\Big(-\frac12\scapro{A^{-1}\vek(O'yO)}{\vek(O' yO)}_{\R^{dJ}}\Big)\,\HH_J(dO)\HH_d(dO'),\; y\in\R^{d\times J}.
\end{align*}
We have the singular-value decomposition $Y=O_Y' L O_Y$ with some orthogonal matrices $O_Y'\in O(d)$, $O_Y\in O(J)$. Then $\Lambda:=LL^\top=\diag(\lambda_1,\ldots,\lambda_d)$ is uniquely determined by $YY^\top$ and
\[ O'YO=(O'O_Y')L(O_YO)\stackrel{d}{=} O'LO\]
holds for independent orthogonal matrices $O\sim \HH_J$, $O'\sim\HH_d$ because of the group invariance of Haar measure. Consequently, by Theorem 3.2.17 and the proof of Theorem 3.2.18 in \cite{muirhead2009} the eigenvalues $\lambda_1\ge\cdots\ge\lambda_d\ge 0$ of $YY^\top$ have density
\begin{align*}
f^\Lambda(\lambda_1,\ldots,\lambda_{d})&=c\det(A)^{-1/2}\prod_{j=1}^{d}\lambda_j^{(J-d-1)/2}\prod_{i<j}(\lambda_i-\lambda_j)\times\\
&
\iint \exp\Big(-\frac12\scapro{A^{-1}\vek(O' L O)}{\vek(O'LO)}_{\R^{dJ}}\Big)\,\HH_J(dO)\HH_d(dO')
\end{align*}
with constant $c=\pi^{d^2/2}2^{-dJ/2}\Gamma(d/2)^{-1}\Gamma(J/2)^{-1}$.
\end{proof}

The preceding density together with a symmetrisation argument yields a stochastic dominance property for the smallest eigenvalue and an explicit bound on its Laplace transform.

\begin{theorem}\label{ThmBernsteinLower}
Let $Y_j\sim N(0,A_j)$, $j=1,\ldots,J$, be independent $d$-dimensional Gaussian random vectors. Suppose that $A_j\ge A_0$ holds for  $j=1,\ldots,J$ with some $A_0\in\R_{spd}^{d\times d}$.
Then we have the stochastic order
\begin{align*}
\forall t>0:\; \PP\Big(\lambda_{min}\Big(\sum_{j=1}^JY_jY_j^\top\Big)\le t\Big) & \le \PP^{(0)}\Big(\lambda_{min}\Big(\sum_{j=1}^JY_jY_j^\top\Big)\le t\Big),
\end{align*}
where $Y_j\sim N(0,A_0)$ i.i.d. holds under $\PP^{(0)}$. Moreover, we can bound the Laplace transform for any $\theta\ge 0$ and $J\ge d$ by
\[ \E\Big[\exp\Big(-\theta \lambda_{min}\Big(\sum_{j=1}^JY_jY_j^\top\Big)\Big)\Big]\le \tfrac{ \Gamma(1/2)\Gamma((J+1)/2)} {\Gamma(d/2)\Gamma((J-d+2)/2)} (1+2\theta\lambda_{min}(A_0))^{-(J-d+1)/2}.
\]
\end{theorem}

\begin{remark}
It seems intuitive that $\lambda_{min}(YY^\top)$ becomes smaller when the $A_j$ are replaced by $A^{(0)}$. This does not follow directly, however, by a standard coupling argument.
\end{remark}

\begin{proof}
Note first that the result is trivial if $A_0$ is not invertible, that is $\lambda_{min}(A_0)=0$. For $J<d$ we always have $\lambda_{min}(\sum_{j=1}^JY_jY_j^\top)=0$ and the stochastic order is immediate. Henceforth, we therefore assume $\lambda_{min}(A_0)>0$ and $J\ge d$.

For our  symmetrisation argument let $F_{\pm 1}=(I_{J}-E_{d,d})\pm E_{d,d}\in \R^{J\times J}$ so that $F_1$ is the identity and $F_{-1}$ flips the sign of the last coordinate. Then by the invariance $O\stackrel{d}{=} F_{\pm 1} O$ for $O\sim\HH_J$  we can introduce a random sign flip $F_\eps$ with a Rademacher random variable $\eps$ in the density formula for $\sum_jY_jY_j^\top$ of Lemma \ref{LemWishartGen} with $Y=(Y_1,\ldots,Y_J)$ to obtain
\begin{align*}
&f^\Lambda(\lambda_1,\ldots,\lambda_{d})=c\prod_{j=1}^J\det(A_j)^{-1/2}\prod_{j=1}^{d}\lambda_j^{(J-d-1)/2}\prod_{i<j}(\lambda_i-\lambda_j)\times\\
&
\quad\iint \E_\eps\Big[\exp\Big(-\frac12\sum_{j=1}^J\scapro{A_j^{-1}O' LF_\eps O_j}{O' LF_\eps O_j} \Big)\Big]\,\HH_J(dO)\HH_d(dO'),
\end{align*}
where $O_j=Oe_j$ is the $j$th column of $O$.
Treating also the matrices $O$ and $O'$ as random (under $\HH_J\otimes\HH_\ell$), we may introduce the joint density
\begin{align*}
f(O,O',\lambda_1,\ldots,\lambda_{d})
:=&c\prod_{j=1}^J\det(A_j)^{-1/2}\prod_{j=1}^{d}\lambda_j^{(J-d -1)/2}\prod_{i<j}(\lambda_i-\lambda_j)\times\\
&\E_\eps\Big[\exp\Big(-\frac12\sum_{j=1}^J\scapro{A_j^{-1}O' LF_\eps O_j}{O' LF_\eps O_j} \Big)\Big].
\end{align*}
In the case $A_j=A_0$ for all $j=1,\ldots,J$ we use $OO^\top=I_{J}=F_\eps^2$ such that the exponent
\begin{align*}
-\frac12\sum_{j=1}^J\scapro{A_0^{-1}O' LF_\eps O_j}{O' LF_\eps O_j} &=-\frac12\trace\Big(A_0^{-1}O' LF_\eps O O^\top F_\eps L^\top (O')^\top\Big)\\
&=
-\frac12\trace\Big(A_0^{-1}O' L L^\top (O')^\top\Big),
\end{align*}
is independent of $\eps\in\{-1,+1\}$. Denoting by $f^{(0)}$ the density $f$ in terms of $A_0$,
the likelihood ratio can be written as
\[\frac{f(O,O',\lambda_1,\ldots,\lambda_{d})}{f^{(0)}(O,O',\lambda_1,\ldots,\lambda_{d})}=\tilde c \E_\eps\Big[ \exp\Big(\frac12\sum_{j=1}^J\scapro{((A^{(0)})^{-1}-A_j^{-1})(O' LF_\eps O_j)}{O' LF_\eps O_j}\Big)\Big]
\]
with some constant $\tilde c>0$. Noting
\[LF_\eps=\Big(\sum_{i=1}^{d-1}\lambda_i^{1/2}E_{i,i}\Big)+\eps \lambda_{d}^{1/2}E_{d,d},\]
the exponent is a quadratic form $Q(\eps\lambda_{d}^{1/2})$ in $\eps \lambda_{d}^{1/2}$. Since $(A^{(0)})^{-1}-A^{-1}$ is positive semi-definite by assumption, we can write $Q(x)=a+bx+cx^2$ with some $a,c\ge 0$ and $b\in\R$. This shows that
\[ \E_\eps[\exp(Q(\eps \lambda_{d}^{1/2}))]= \exp(a+c\lambda_{d})\cosh(b\lambda_{d}^{1/2})\]
is always increasing in $\lambda_{d}\ge 0$. Hence,  the likelihood ratio $f/f^{(0)}$ is increasing in $\lambda_{d}$. Integrating $O,O'$ and $\lambda_i$ for $i\le d-1$ out, we thus deduce that $\lambda_d(\sum_{j=1}^JY_jY_j^\top)$ is stochastically larger under $\PP$ than under $\PP^{(0)}$, which is the stochastic order result.

In view of $A_j\ge \lambda_{min}(A_0)I_d$ the stochastic order yields further for the Laplace transform with $\zeta_j\sim N(0,I_d)$ i.i.d.
\[ \E\Big[\exp\Big(-\theta\lambda_{min}\Big(\sum_{j=1}^{J}Y_jY_j^\top\Big)\Big)\Big]\le \E\Big[\exp\Big(-\theta\lambda_{min}(A_0)\lambda_{min}\Big(\sum_{j=1}^{J}\zeta_j\zeta_j^\top\Big)\Big)\Big],\quad \theta\ge 0.
\]
Using the bound for the density of $\lambda_{min}$ in the proof of Proposition 5.1 in \citet{Edelman1988}, we obtain for $\tilde\theta\ge 0$
\begin{align*}
\E\Big[&\exp\Big(-\tilde\theta\lambda_{min}\Big(\sum_{j=1}^{J}\zeta_j\zeta_j^\top\Big)\Big)\Big]\\
&\le \tfrac{ 2^{-(J-d+1)/2} \Gamma(1/2)\Gamma((J+1)/2)} {\Gamma(d/2)\Gamma((J-d+1)/2)\Gamma((J-d+2)/2)} \int_0^\infty \lambda^{(J-d-1)/2}e^{-\lambda(1/2+\tilde\theta)}d\lambda\\
&= \tfrac{ \Gamma(1/2)\Gamma((J+1)/2)} {\Gamma(d/2)\Gamma((J-d+2)/2)} (1+2\tilde\theta)^{-(J-d+1)/2}.
\end{align*}
With $\tilde\theta=\theta\lambda_{min}(A_0)$ this gives the asserted bound.
\end{proof}

For a sufficiently large number of summands the previous bound can be generalised to the setting where each summand has a population covariance where only the size of the eigenvalue is lower bounded.

\begin{corollary}\label{CorEigenvalueLB}
Let $Y_j\sim N(0,A_j)$, $j=1,\ldots,J$, be independent $d$-dimensional Gaussian random vectors with $\lambda_\ell(A_j)\ge\underline\lambda_\ell>0$ for all $j$ and some $\ell\in\{1,\ldots,d\}$. Set $J_0:=\ceil{J/C_{d,\ell}}$ with a fixed constant $C_{d,\ell}>1$ only depending on $\ell$ and $d$  and assume $J_0\ge \ell$. Then  we have
the Laplace transform bound for $\theta\ge 0$:
\[ \E\Big[\exp\Big(-\theta \lambda_{\ell}\Big(\sum_{j=1}^JY_jY_j^\top\Big)\Big)\Big]\le \tfrac{ \Gamma(1/2)\Gamma(J_0+1)/2)} {\Gamma(\ell/2)\Gamma((J_0-\ell+2)/2)} (1+\theta\underline\lambda_\ell/2)^{-(J_0-\ell+1)/2}.
\]
\end{corollary}

\begin{proof}
Denote by $P_{j,\le \ell}$ the orthogonal projection onto the $\ell$-dimensional eigenspace of $A_j$ corresponding to $\lambda_1(A_j),\ldots,\lambda_\ell(A_j)$. By the metric entropy result for Grassmannian manifolds \citep[Proposition 6]{Pajor1998} and the relationship with internal covering numbers there is a family $\mathcal V$ of $\ell$-dimensional subspaces of $\R^d$ such that with orthogonal projections $P_V$ onto $V$
\[ \forall j=1,\ldots,J\,\exists V\in{\mathcal V}:\; \norm{P_V-P_{j,\le\ell}}\le 1/2\]
and $\mathcal V$ has cardinality less than $C_{d,\ell}:=C^{\ell(d-\ell)}$ for some universal constant $C> 1$. Hence, by a counting argument there are $V\in{\mathcal V}$ and ${\mathcal J}_0\subset\{1,\ldots,J\}$ with $\norm{P_V-P_{j,\le\ell}}\le 1/2$ for all $j\in{\mathcal J}_0$ with $\abs{{\mathcal J}_0}= \ceil{JC_{d,\ell}^{-1}}=J_0$. For $v\in V$ and $j\in{\mathcal J}_0$ we then obtain
\begin{align*}
\scapro{A_jv}{v} &\ge \lambda_\ell(A_j)\scapro{P_{j,\le\ell}v}{v}\ge \underline\lambda_\ell \norm{v+(P_{j,\le\ell}-P_V)v}^2\\
& \ge \underline\lambda_\ell \big(\norm{v}-\norm{(P_{j,\le\ell}-P_V)v}\big)^2\ge \tfrac{\underline\lambda_\ell}4\norm{v}^2.
\end{align*}
We apply Theorem \ref{ThmBernsteinLower} in dimension $\ell$ for the restrictions $P_VY_j$, $P_VA_j|_{V}$ to the subspace $V$, for $j\in{\mathcal J}_0$, $\abs{{\mathcal J}_0}=J_0\ge\ell$, and with $A_0=\tfrac{\underline\lambda_\ell}4 \Id_V$ and obtain
\begin{align*}
\E\Big[\exp\Big(-\theta \lambda_{\ell}\Big(\sum_{j=1}^JY_jY_j^\top\Big)\Big)\Big] &\le \E\Big[\exp\Big(-\theta \lambda_{min}\Big(\sum_{j\in{\mathcal J}_0}P_VY_jY_j^\top|_{V}\Big)\Big)\Big]\\
&\le \tfrac{ \Gamma(1/2)\Gamma(J_0+1)/2)} {\Gamma(\ell/2)\Gamma((J_0-\ell+2)/2)} (1+\theta\underline\lambda_\ell/2)^{-(J_0-\ell+1)/2},
\end{align*}
as claimed.
\end{proof}

Given the Laplace transform result we obtain a deviation inequality for a triangular scheme over different Wishart matrices. This is exactly what we need for analysing our test statistics $T_{n,h}$ under the alternatives.

\begin{corollary}\label{CorEigenvalueTriangular}
Let $Y_{jk}\sim N(0,A_{jk})$, $j=1,\ldots,J$, $k=1,\ldots,K$, be independent $d$-dimensional random vectors  with $\underline\lambda_{\ell,k}:=\min_j\lambda_\ell(A_{jk})>0$ for all $k$ and some $\ell\in\{1,\ldots,d\}$. Without loss of generality suppose the order $\underline\lambda_{\ell,1}\ge \underline\lambda_{\ell,2}\ge\cdots\ge\underline\lambda_{\ell,K}$. Set $J_0:=\ceil{J/C_{d,\ell}}$ and assume $J_0\ge 2\ell$ with the constant $C_{d,\ell}$ from Corollary \ref{CorEigenvalueLB}. Then  we have for all $\tau>0$ and $K'=1,\ldots,K$ the lower tail bound
\begin{align*}
\PP\Big(\sum_{k=1}^K \lambda_\ell\Big(\frac1{J}\sum_{j=1}^JY_{jk}Y_{jk}^\top\Big)\le \tau K'\underline\lambda_{\ell,K'}\Big)\le 
\Big(\bar C_{d,\ell} \tau\Big)^{K'(J_0-\ell+1)/2},
\end{align*}
where the constant $\bar C_{d,\ell}>1$ only depends on $d$ and $\ell$.
\end{corollary}

\begin{remark}
Let us note that $\max_{K'}K'\underline\lambda_{\ell,K'}$ is the weak-$\ell^1$ norm of the vector $(\underline\lambda_{\ell,k})_k$.
\end{remark}

\begin{proof}
By elementary properties of binomial coefficients and Gamma functions we have $\bar C:=\sup_{J_0\ge2\ell}(\frac{ \Gamma(1/2)\Gamma((J_0+1)/2)} {\Gamma(\ell/2)\Gamma((J_0-\ell+2)/2)})^{2/(J_0-\ell+1)}<\infty$. We thus infer from Corollary \ref{CorEigenvalueLB}
\begin{align*}
\E\Big[\exp\Big(-\theta \lambda_{\ell}\Big(\sum_{j=1}^JY_{jk}Y_{jk}^\top\Big)\Big)\Big]
&\le \Big(\bar C^{-1}\theta \underline\lambda_{\ell,k}/2
\Big)^{-(J_0-\ell+1)/2},\quad \theta\ge 0.
\end{align*}
The classical Chernoff argument yields for $t>0$ and any $K'\le K$
\begin{align*}
\PP\Big(\sum_{k=1}^K \lambda_\ell\Big(\frac1{J}\sum_{j=1}^JY_{jk}Y_{jk}^\top\Big)\le \tau K'\underline\lambda_{\ell,K'}\Big)
&\le e^{\theta \tau K'\underline\lambda_{\ell,K'}  }\prod_{k=1}^{K}\E\Big[\exp\Big(-\frac \theta J \lambda_\ell\Big(\sum_{j=1}^JY_{jk}Y_{jk}^\top\Big)\Big)\Big]\\
&\le e^{\theta \tau K'\underline\lambda_{\ell,K'}}\prod_{k=1}^{K'} \Big(\bar C^{-1}\theta J^{-1}\underline\lambda_{\ell,k}/2\Big)^{-(J_0-\ell+1)/2}\\
&\le e^{\theta \tau K'\underline\lambda_{\ell,K'}}\Big(\bar C^{-1}\theta J^{-1}\underline\lambda_{\ell,K'}/2\Big)^{-(J_0-\ell+1)K'/2}.
\end{align*}
Choosing $\theta=(J_0-\ell+1)/(2\tau\underline\lambda_{\ell,K'})$ and applying $J\le C_{d,\ell}J_0$, $J_0-\ell+1\ge J_0/\ell$ gives the asserted bound with $\bar C_{d,\ell}=4e\ell \bar C C_{d,\ell}$.
\end{proof}

\subsection{Technical results for Section \ref{SecResults}} \label{AppTechnical3}

\begin{proposition}\label{PropCalcCritVal}
Consider the test statistics  $T_{n,h}$ in \eqref{Eqphi0} and assume $\rank(\Sigma_X(t))\le r$ for all $t$. Let
\[\Delta_p(\Sigma_X,h):=\Big(\sum_{k=0}^{h^{-1}-1}h\Delta_p(\Sigma_X,I_k)^p\Big)^{1/p}\]
denote the average $L^p$-variation of $\Sigma_X$ over blocks $I_k$ and introduce the constants
\[C_{\rho,1}=(1+2\sqrt{\log \rho} ) \sqrt{2\rho+4},\quad C_{\rho,2}=(1+\log\rho)(\rho+4),\quad\rho\in\N.\]
Then we have with probability at least $1-\alpha\in(0,1)$ for any $\delta>0$ and $p\ge 2$:
\begin{enumerate}
\item without a spectral gap assumption
\begin{align*}
T_{n,h}
& \le (1+\delta)\Big( \big(\Delta_1(\Sigma_X,h)+\norm{\Sigma_Z}_{L^1}\big)
  +C_{d-r,1}(2\Delta_2(\Sigma_X,h)+\norm{\Sigma_{Z}}_{L^2})(nh)^{-1/2}\\
 &\quad +4
 \big(2\Delta_p(\Sigma_X,h)+\norm{\Sigma_Z}_{L^p}\big)\big(C_{d-r,2}(nh)^{-1+1/p}+\delta^{-1}n^{-1+1/p}\log(\alpha^{-1})\big)\Big);
\end{align*}

\item under the spectral gap assumption $\underline\lambda_{r}=\inf_{t\in[0,1]}\lambda_r(\Sigma_X(t))>0$
\begin{align*}
T_{n,h}
 &\le (1+\delta)\Big( \big(\tfrac{2}{\underline\lambda_r}\Delta_2(\Sigma_{X},h)^2+\norm{\Sigma_{Z}}_{L^1}\big)\\
 &\quad
  +C_{d-r,1} \big(\tfrac{r+4}{\underline\lambda_r }  \Delta_4(\Sigma_X,h)^{2} + \norm{\Sigma_Z}_{L^2}\big)(nh)^{-1/2}\\
 &\quad +4\big(\tfrac{r+4}{\underline\lambda_r }\Delta_{2p}(\Sigma_X,h)^2+\norm{\Sigma_Z}_{L^{p}}\big)\big(C_{d-r,2}(nh)^{-1+1/p}+\delta^{-1}n^{-1+1/p}\log(\alpha^{-1})\big) \Big).
\end{align*}

\end{enumerate}
\end{proposition}

\begin{remark}
It is implicitly understood that the appearing $L^p$-variations and $L^p$-norms of $\Sigma_X$ and $\Sigma_Z$ are all  finite.

The bounds consist in each case of four different contributions. The first and asymptotically dominant term comes from the deterministic bias. The second and third  capture the additional bias induced by the expected $(r+1)$st eigenvalue of $\hat\Sigma^{kh}_Y$ and scales in the number of observations on a block like $(nh)^{-1/2}$ for the subgaussian deviations and $(nh)^{-1+1/p}$ for the subexponential deviations. The dependence on $p$ comes from bounding the maximum of blockwise integrals by an $L^p$-norm, similarly to Sobolev embeddings, and one might think of $p=\infty$ for a first intuition. The size of random fluctuations scales with the total number of observations so that the fourth term is  of order  $n^{-1+1/p}$. In the case of a spectral gap, the squared error structure induces a natural $L^{2p}$-variation bound. The proof reflects exactly this error decomposition. Concerning the dimension dependence we conjecture that the linearity up to logarithmic terms of $C_{d-r,1}^2$ and $C_{d-r,2}$ in the remaining dimension $d-r$ is also necessary.
\end{remark}

\begin{proof}\mbox{}
To establish part (a), let us consider for each block $k$ the  eigenspace $V_{>r}:=\spann(v_{r+1},\ldots,v_d)$ of eigenvectors corresponding to the $d-r$ smallest eigenvalues  $\lambda_{r+1}( \Sigma_X^{kh})\ge\cdots\ge\lambda_d( \Sigma_X^{kh})\ge 0$ of $\Sigma_X^{kh}$. With the orthogonal projection $P_{>r}$ onto $V_{>r}$ we introduce the $(d-r)\times (d-r)$-lower right minors
\begin{equation}\label{EqMinor}
S_{>r}:=P_{>r}S|_{V_{>r}}\text{ for matrices } S\in\{\hat \Sigma_X^{kh},\Sigma_X^{kh},\hat \Sigma_Y^{kh},\Sigma_Y^{kh}\}.
\end{equation}
By the Cauchy interlacing law (e.g., \citet{johnstone2001} or \citet{tao2012}),
\begin{align*}
\lambda_{r+1}(\hat\Sigma_Y^{kh})&\le \lambda_{max}(\hat \Sigma^{kh}_{Y,>r})\\
&\le \lambda_{max}(\hat \Sigma^{kh}_{Y,>r}-\Sigma^{kh}_{Y,>r})+\lambda_{max}(\Sigma_{Y,>r}^{kh})\\
&\le \lambda_{max}(\hat \Sigma^{kh}_{Y,>r}-\Sigma^{kh}_{Y,>r})+\norm{\Sigma^{kh}_{Z}}+\lambda_{r+1}(\Sigma_X^{kh}).
\end{align*}

Introduce $I_{jk}=[kh+(j-1)/n,kh+j/n]$. We apply the matrix deviation inequality for triangular schemes from Theorem \ref{ThmBernsteinTriangular} in dimension $d-r$ with $A_{jk}=\int_{I_{jk}}\Sigma_{Y,>r}(t)dt$, $A_k=h\Sigma_{Y,>r}^{kh}$. From Corollary \ref{CorBias0} we obtain
\[\norm{A_k}\le h\big(\norm{\Sigma_{X,>r}^{kh}}+\norm{\Sigma_{Y,>r}^{kh}-\Sigma_{X,>r}^{kh}}\big)\le h\Delta_1(\Sigma_{X},I_k)+\norm{\Sigma_{Z}}_{L^1(I_k)}
\]
 in the absence of a spectral gap.  Then Proposition \ref{PropNormPerturb2} and   H\"older's inequality yield
\begin{align*}
\norm{A_{jk}} &\le \int_{I_{jk}} \big(\norm{\Sigma_{X,>r}(t)}+\norm{\Sigma_{Z}(t)}\big)\,dt \\
&\le n^{-1}\Big(2(nh)^{1/p}\Delta_p(\Sigma_X,I_{k})+ n^{1/p}\norm{\Sigma_{Z}}_{L^p(I_{jk})}\Big)
\end{align*}
and by Jensen's inequality
\begin{align*}
\sum_{j=1}^{nh} \norm{A_{jk}}^2 &\le \sum_{j=1}^{nh} \frac1n\int_{I_{jk}}\big(\norm{\Sigma_{X,>r}(t)}+ \norm{\Sigma_{Z}(t)}\big)^2 dt\\
&\le \frac1n\int_{I_k} \big(\Delta_1(\Sigma_X,I_{k})+\norm{\Sigma_X(t)-\Sigma_X^{kh}} +\norm{\Sigma_{Z}(t)}\big)^2dt\\
&\le \frac {2}n\Big(3h \Delta_1(\Sigma_X,I_{k})^2+ h \Delta_2(\Sigma_X,I_{k})^2+\norm{\Sigma_{Z}}_{L^2(I_k)}^2 \Big)\\
&\le 2n^{-1}\Big(4h\Delta_2(\Sigma_X,I_k)^2+\norm{\Sigma_{Z}}_{L^2(I_k)}^2\Big).
\end{align*}
The quantities in Theorem \ref{ThmBernsteinTriangular} are therefore bounded as
\begin{align*}
R_k &\le (d-r+4)n^{-1+1/p}\Big(2h^{1/p}\Delta_p(\Sigma_X,I_{k})+\norm{\Sigma_{Z}}_{L^p(I_k)}\Big),\\
\sigma_k^2 &\le 2(d-r+2)n^{-1}\Big(4h\Delta_2(\Sigma_X,I_k)^2+\norm{\Sigma_{Z}}_{L^2(I_k)}^2\Big).
\end{align*}
By Theorem \ref{ThmBernsteinTriangular} with $t=\log(\alpha^{-1})$  we thus have with probability at least $1-\alpha$ for any $\delta>0$
\begin{align*}
T_{n,h}
&\le (1+\delta)\Big(\big(\Delta_1(\Sigma_X,h)+\norm{\Sigma_Z}_{L^1}\big)+(2\sqrt{\log (d-r)}+1)\Big(h^{-1}\sum_{k=0}^{h^{-1}-1}\sigma_k^2\Big)^{1/2}\\
&\quad +4(\log(d-r)+1)\sum_{k=0}^{h^{-1}-1}R_k\Big) + 4(1+\delta^{-1})n^{-1+1/p}\big(2\Delta_p(\Sigma_X,h)+\norm{\Sigma_Z}_{L^p}\big) \log(\alpha^{-1})\\
 &\le (1+\delta) \Big(\big(\Delta_1(\Sigma_X,h)+\norm{\Sigma_Z}_{L^1}\big)\\
 &\quad +(2\sqrt{\log (d-r)} +1) \Big(2(d-r+2)(nh)^{-1}(4\Delta_2(\Sigma_X,h)^2+\norm{\Sigma_{Z}}_{L^2}^2)\Big)^{1/2}\\
 &\quad +4(\log(d-r)+1)(d-r+4)(nh)^{-1+1/p}\big(2\Delta_p(\Sigma_X,h)+\norm{\Sigma_Z}_{L^p}\big)\\
&\quad + 4\delta^{-1}n^{-1+1/p}\big(2\Delta_p(\Sigma_X,h)+\norm{\Sigma_Z}_{L^p}\big) \log(\alpha^{-1})\Big)\\
 &\le (1+\delta)\Big( \big(\Delta_1(\Sigma_X,h)+\norm{\Sigma_Z}_{L^1}\big)+(nh)^{-1/2}C_{d-r,1}(2\Delta_2(\Sigma_X,h)+\norm{\Sigma_{Z}}_{L^2})\\
 &\quad +4n^{-1+1/p}\big(C_{d-r,2}h^{-1+1/p}+\delta^{-1}\log(\alpha^{-1})\big)
 \big(2\Delta_p(\Sigma_X,h)+\norm{\Sigma_Z}_{L^p}\big)\Big).
\end{align*}

In the case (b) of a spectral gap we use $A_{jk}$ and $A_k$ as before, but employ the quadratic matrix deviation bounds. Then Corollary \ref{CorBias0} yields
    \[\norm{A_k}\le h\big(\norm{\Sigma_{X,>r}^{kh}}+\norm{\Sigma_{Y,>r}^{kh}-\Sigma_{X,>r}^{kh}}\big)\le \tfrac{2} {\underline\lambda_r}h\Delta_2(\Sigma_{X},I_k)^2+\norm{\Sigma_{Z}}_{L^1(I_k)}.
\]
By  Proposition \ref{PropNormPerturb2} and H\"older's inequality we have
\[ \norm{A_{jk}} \le n^{-1}\Big( \frac{r+4}{\underline\lambda_r } (nh)^{1/p}\Delta_{2p}(\Sigma_X,I_k)^{2} + n^{1/p}\norm{\Sigma_Z}_{L^{p}(I_{jk})} \Big)
\]
as well as
\[ \sum_{j=1}^{nh}\norm{A_{jk}}^2 \le 2n^{-1}\Big(\frac{(r+4)^2}{\underline\lambda_r^2 }  h\Delta_4(\Sigma_X,I_k)^{4} + \norm{\Sigma_Z}_{L^2(I_k)}^2\Big).
\]
The quantities in Theorem \ref{ThmBernsteinTriangular} are in this case bounded as
\begin{align*}
R_k &\le (d-r+4)n^{-1+1/p}\Big( \frac{r+4}{\underline\lambda_r } h^{1/p}\Delta_{2p}(\Sigma_X,I_k)^{2} + \norm{\Sigma_Z}_{L^{p}(I_{k})} \Big),\\
\sigma_k^2 &\le 2(d-r+2)n^{-1}\Big(\frac{(r+4)^2}{\underline\lambda_r^2 }  h\Delta_4(\Sigma_X,I_k)^{4} + \norm{\Sigma_Z}_{L^2(I_k)}^2\Big).
\end{align*}
We apply Theorem \ref{ThmBernsteinTriangular} exactly as before and arrive at the result in (b).
\end{proof}

\begin{proof}[Proof of Proposition \ref{PropLB}.]
Let us first treat the case $d=2$, $r=1$, $L=4\pi$ and $\underline\lambda_1\ge n^{-\beta}/\sqrt 2$.
Consider the covariance function $\Sigma_X(t)$ from Example \ref{ExRank} for $h=n^{-1}$, satisfying  $\Sigma_X\in {\mathcal H}_0^{gap}(1,\beta,4\pi,0,\underline\lambda_{1})$. On the other hand, the constant covariance $\tilde\Sigma_X(t):=\diag(\underline\lambda_1,(2\underline\lambda_1)^{-1}n^{-2\beta})$ lies in the alternative ${\mathcal H}_1(1,\hbar,(2\underline\lambda_1)^{-1}n^{-2\beta})$ even for $\hbar=1$. Because of $\int_{(i-1)/n}^{i/n}\Sigma_X(t)dt=\int_{(i-1)/n}^{i/n}\tilde\Sigma_X(t)dt$ for all $i$, the observations $X(i/n)$, $i=0,\ldots,n$, have the same law under $\Sigma_X$ and $\tilde\Sigma_X$ which entails  $\E_{\Sigma_X}[\phi]=\E_{\tilde\Sigma_X}[\phi]$  for any test $\phi$, as asserted.
In the case $\underline\lambda_1\in (0, n^{-\beta}/\sqrt 2)$ observe the inclusions
\[ {\mathcal H}_0^{gap}(1,\beta,L,0,n^{-\beta}/\sqrt 2)\subset {\mathcal H}_0^{gap}(1,\beta,L,0,\underline\lambda_{1}) \subset {\mathcal H}_0(1,\beta,L,0).\]
Therefore, the result for $\underline\lambda_1=n^{-\beta}/\sqrt 2$, where $\tilde\Sigma_X\in {\mathcal H}_1(1,\hbar,2^{-1/2}n^{-\beta})$, extends to this case. Altogether we  have thus shown  the assertion for $d=2$, $r=1$, $L=4\pi$.

We reduce the general case to this particular choice. Indeed, for general $1\le r<d$ it suffices to consider covariance functions $\Sigma_X(t)\in\R^{d\times d}$ that are diagonal with a large constant (larger than $\underline\lambda_r$) in the coordinates $\ell=1,\ldots,r-1$ and equal to zero in the coordinates $\ell=r+2,\ldots,d$, while equal to the above $(2\times2)$-covariance functions in the entries with coordinates $\ell\in\{r,r+1\}$. Then all H\"older and eigenvalue conditions are clearly fulfilled. Moreover, given that the result is proved for $\Sigma_X$ with the H\"older constant $L=4\pi$, then the covariance $\frac{L}{4\pi}\Sigma_X(t)$ with $\Sigma_X(t)$ from above is in $C^\beta(L)$. The eigenvalue bound $\underline\lambda_r$  scales with the factor $\frac{L}{4\pi}$ as well. The same holds for $\frac{L}{4\pi}\tilde\Sigma_X$ and the general assertion follows by a simple rescaling argument.
\end{proof}

\subsection{Technical results for Section \ref{SecStochVol}}\label{AppStochVol}

Throughout we work under the model \eqref{EqStochVol} and Assumption \ref{AssStochVol} for $p\in\{1,2\}$ and with the blockwise estimators $\hat\Sigma_X^{kh'}$ where $nh',(h')^{-1}\in\N$. We first show two approximation results before establishing the Central Limit Theorem \ref{ThmNQVCLT}.

\begin{lemma}\label{LemVolVol}
Uniformly over $k\in\{1,\ldots,(h')^{-1}-2\}$
\begin{align*}
\bnorm{\hat\Sigma_X^{(k+1)h'}-2\hat\Sigma_X^{kh'}+\hat\Sigma_X^{(k-1)h'}}
&=\bnorm{\Gamma((k-1)h')\int_{(k-1)h'}^{(k+2)h'}w_k(s)dB'(s)}\\
&\quad +{\mathcal O}_{L^{2p}}\Big((h')^{(\beta_b+1)\wedge (\beta_\Gamma+1/2)}+(nh')^{-1/2}\Big),
\end{align*}
where for $s\in[(k-1)h',(k+2)h']$
\begin{align*}
w_k(s)&=\frac1{h'}\int_{0}^{h'}\Big({\bf 1}\big(s-u\in[kh',(k+1)h']\big)-{\bf 1}\big(s-u\in[(k-1)h',kh']\big)\Big)\,du.
\end{align*}
\end{lemma}

\begin{proof}
From Equation \eqref{EqCovSigmakh}, conditioning on $\Sigma_X$, we infer
\[ \E\Big[\norm{\hat\Sigma_X^{kh'}-\Sigma_X^{kh'}}^{2p}\Big]^{1/2p}\lesssim (nh')^{-1/2},\]
using $\Sigma_0,b(t),\Gamma(t)\in L^{2p}$ and thus $\max_{t\in[0,1]}\E[\norm{\Sigma_X(t)}^{2p}]\lesssim 1$. Applying this to $\hat\Sigma_X^{(k+1)h'}$ as well, we can  expand
\begin{align*}
\hat\Sigma_X^{(k+1)h'}-\hat\Sigma_X^{kh'}&=
\frac1{h'}\int_0^{h'}\int_{kh'}^{(k+1)h'}\big(b(t+u)dt+\Gamma(t+u)dB'(t+u)\big)du+{\mathcal O}_{L^{2p}}((nh')^{-1/2}).
\end{align*}
Together with the expansion for $\hat\Sigma_X^{kh'}-\hat\Sigma_X^{(k-1)h'}$ and the regularity condition on $b$ we obtain
\begin{align*}
\hat\Sigma_X^{(k+1)h'}-2\hat\Sigma_X^{kh'}+\hat\Sigma_X^{(k-1)h'}
&=\int_{(k-1)h'}^{(k+2)h'}w_k(s)\Gamma(s)dB'(s)+{\mathcal O}_{L^{2p}}\Big((h')^{1+\beta_b}+(nh')^{-1/2}\Big).
\end{align*}
The regularity of $\Gamma$ yields via the Burkholder-Davis-Gundy inequality
\[ \int_{(k-1)h'}^{(k+2)h'}w_k(s)\Gamma(s)dB'(s)=\Gamma((k-1)h')\int_{(k-1)h'}^{(k+2)h'}w_k(s)dB'(s)+{\mathcal O}_{L^{2p}}((h')^{\beta_\Gamma+1/2}).
\]
All constants depend uniformly in $k$ on the  regularity conditions in Assumption \ref{AssStochVol}.
\end{proof}

Lemma \ref{LemVolVol} allows to derive the asymptotics of the conditional expectations and variances of the normed second differences.

\begin{proposition}\label{PropCondExpSigmaHat}
We have
\begin{align*} &\sum_{k=0}^{(3h')^{-1}-1}\E\Big[\frac{3h'}{(2h')^{p/2}}\norm{\hat\Sigma_X^{(3k+2)h'}-2\hat\Sigma_X^{(3k+1)h'}+\hat\Sigma_X^{3kh'}}^p\,\Big|\,{\mathcal F}_{3kh'}\Big]\\
 &\quad =\int_0^1 \rho_p(\Gamma(t))\,dt +{\mathcal O}_{L^2}\Big((h')^{(\beta_b+1/2)\wedge \beta_\Gamma}+(h')^{-1}n^{-1/2}\Big)
\end{align*}
as well as
\begin{align*} &\sum_{k=0}^{(3h')^{-1}-1}\Var\Big(\frac{3h'}{(2h')^{p/2}}\norm{\hat\Sigma_X^{(3k+2)h'}-2\hat\Sigma_X^{(3k+1)h'}+\hat\Sigma_X^{3kh'}}^p\,\Big|\,{\mathcal F}_{3kh'}\Big)\\
 &\quad =h'\Big(
 3\int_0^1 \big(\rho_{2p}(\Gamma(t))-\rho_p(\Gamma(t))^2\big)\,dt+{\mathcal O}_{L^1}\Big((h')^{(\beta_b+1/2)\wedge \beta_\Gamma}+(h')^{-1}n^{-1/2}\Big)\Big).
\end{align*}
\end{proposition}

\begin{proof}
We inject the value $\norm{w_k}_{L^2}^2=2h'$ into Lemma \ref{LemVolVol} and obtain
\begin{align*}
&\sum_{k=0}^{(3h')^{-1}-1}(h')^{1-p/2}\E\Big[\norm{\hat\Sigma_X^{(3k+2)h'}-2\hat\Sigma_X^{(3k+1)h'}+\hat\Sigma_X^{3kh'}}^p\,\Big|\,{\mathcal F}_{3kh'}\Big]\\
&=\sum_{k=0}^{(3h')^{-1}-1}(h')^{1-p/2}(2h')^{p/2}\rho_p(\Gamma(3kh'))+{\mathcal O}_{L^2}\Big((h')^{(\beta_b+1/2)\wedge \beta_\Gamma}+(h')^{-1}n^{-1/2}\Big).
\end{align*}
This is immediate for $p=1$ and follows for $p=2$ by $\Gamma((k-1)h')\int_{(k-1)h'}^{(k+2)h'}w_k(s)dB'(s)={\mathcal O}_{L^{2p}}((h')^{1/2})$ and an application of Cauchy-Schwarz inequality to the squared norm in Lemma \ref{LemVolVol}.
In the proof of Proposition \ref{PropNQV} $\abs{\rho_p(\Gamma(t))-\rho_p(\Gamma(s))}={\mathcal O}_{L^2}(\abs{t-s}^{\beta_\Gamma})$ was established so that a Riemann sum approximation yields
\[ \sum_{k=0}^{(3h')^{-1}-1}3h'\rho_p(\Gamma(3kh'))=\int_0^1 \rho_p(\Gamma(t))\,dt+{\mathcal O}_{L^2}((h')^{\beta_\Gamma}).
\]
Putting the asymptotics together yields the conditional expectation result.

Arguing similarly for the conditional variances,  we calculate
\begin{align}
&\sum_{k=0}^{(3h')^{-1}-1}(h')^{1-p}\E\Big[\norm{\hat\Sigma_X^{(3k+2)h'}-2\hat\Sigma_X^{(3k+1)h'}+\hat\Sigma_X^{3kh'}}^{2p}\,\Big|\,{\mathcal F}_{3kh'}\Big]\nonumber\\
&=\sum_{k=0}^{(3h')^{-1}-1}(h')^{1-p}(2h')^p\rho_{2p}(\Gamma(3kh'))+{\mathcal O}_{L^1}\Big((h')^{(\beta_b+1/2)\wedge \beta_\Gamma}+(h')^{-1}n^{-1/2}\Big)\nonumber\\
&= \frac{2^p}3\int_0^1 \rho_{2p}(\Gamma(t))\,dt+{\mathcal O}_{L^1}\Big((h')^{(\beta_b+1/2)\wedge \beta_\Gamma}+(h')^{-1}n^{-1/2}\Big).\label{EqNorm2p}
\end{align}
On the other hand, we obtain for the squared conditional expectations
\begin{align*}
&\sum_{k=0}^{(3h')^{-1}-1}(h')^{1-p}\E\Big[\norm{\hat\Sigma_X^{(3k+2)h'}-2\hat\Sigma_X^{(3k+1)h'}+\hat\Sigma_X^{3kh'}}^p\,\Big|\,{\mathcal F}_{3kh'}\Big]^2\\
&=\sum_{k=0}^{(3h')^{-1}-1}(h')^{1-p}(2h')^p\rho_p(\Gamma(3kh'))^2+{\mathcal O}_{L^1}\Big((h')^{(\beta_b+1/2)\wedge \beta_\Gamma}+(h')^{-1}n^{-1/2}\Big)\\
&= \frac{2^p}3\int_0^1 \rho_p(\Gamma(t))^2\,dt+{\mathcal O}_{L^1}\Big((h')^{(\beta_b+1/2)\wedge \beta_\Gamma}+(h')^{-1}n^{-1/2}\Big).
\end{align*}
The difference of the expressions gives the result for the variance.
\end{proof}

\begin{proof}[Proof of Theorem \ref{ThmNQVCLT}]
First observe that $h'n^{1/3}\to\infty$ implies $(h')^{-1}n^{-1/2}=o((h')^{1/2})$ and the ${\mathcal O}$-terms in Proposition \ref{PropCondExpSigmaHat} are  asymptotically negligible.

We use the conditional expectation and variance results from Proposition \ref{PropNQV} and the fact that $\norm{\hat\Sigma_X^{(3k+2)h'}-2\hat\Sigma_X^{(3k+1)h'}+\hat\Sigma_X^{3kh'}}^2$ is an even function of Brownian increments to apply a standard stable limit theorem (Theorem 4.2.1 in \citet{JacPro2011}, compare also its application there in Theorem 5.3.5(i-$\alpha$)),  noting that we need not require the differentiability of the functional because the linear term becomes already negligible by our choice of $h'$.
\end{proof}

\begin{proof}[Proof of Proposition \ref{PropVarhat}]
By Lemma \ref{LemVolVol} and $\norm{\hat\Sigma_X^{(3k+2)h'}-2\hat\Sigma_X^{(3k+1)h'}+\hat\Sigma_X^{3kh'}}={\mathcal O}_{L^{2p}}((h')^{1/2})$, uniformly in $k$, we find
\begin{align*}
&\sum_{k=0}^{(3h')^{-1}-1} (h')^{1-p} \norm{\hat\Sigma_X^{(3k+2)h'}-2\hat\Sigma_X^{(3k+1)h'}+\hat\Sigma_X^{3kh'}}^{2p}\\
&\quad= \Big(\sum_{k=0}^{(3h')^{-1}-1} (h')^{1-p}\bnorm{\Gamma(3kh')\int_{3kh'}^{(3k+3)h'}w_{3k+1}(s)dB'(s)}^{2p}\Big)\\
&\qquad+{\mathcal O}_{L^1}\Big((h')^{(\beta_b+1/2)\wedge \beta_\Gamma}+(h')^{-1}n^{-1/2}\Big)
\end{align*}
From \eqref{EqNorm2p} we  obtain
\begin{align*} &\sum_{k=0}^{(3h')^{-1}-1}  \E\Big[\frac{3h'}{(2h')^p}\bnorm{\Gamma(3kh')\int_{3kh'}^{(3k+3)h'}w_{3k+1}(s)dB'(s)}^{2p}\,\Big|\,{\mathcal F}_{3kh'}\Big]\\
&\quad= \int_0^1\rho_{2p}(\Gamma(t))dt+{\mathcal O}_{L^1}\Big((h')^{1+\beta_\Gamma}\Big).
\end{align*}
Looking at $4p$-moments, we deduce directly that
\begin{align*}
&\sum_{k=0}^{(3h')^{-1}-1}\Var\Big(\frac{3h'}{(2h')^p} \bnorm{\Gamma(3kh')\int_{3kh'}^{(3k+3)h'}w_{3k+1}(s)dB'(s)}^{2p}\,\Big|{\mathcal F}_{3kh'}\Big)\\
&\quad \le  9\sum_{k=0}^{(3h')^{-1}-1}\rho_{4p}(\Gamma(3kh'))(h')^2,
\end{align*}
which is ${\mathcal O}_P(h')$ since $\max_{t\in[0,1]}\E[\norm{\Gamma(t)}^{4p}]<\infty$. Standard martingale arguments \cite[Lemma 2.2.11(a)]{JacPro2011} therefore yield
\[\frac{3h'}{(2h')^p}\sum_{k=0}^{(3h')^{-1}-1} \norm{\hat\Sigma_X^{(3k+2)h'}-2\hat\Sigma_X^{(3k+1)h'}+\hat\Sigma_X^{3kh'}}^{2p} \xrightarrow{\PP} \int_0^1\rho_{2p}(\Gamma(t))dt
\]
and thus  the first limiting result.

For the second result the same arguments give
\begin{align*}
&\sum_{k=0}^{(6h')^{-1}-1} (h')^{1-p}\norm{\hat\Sigma_X^{(6k+2)h'}-2\hat\Sigma_X^{(6k+1)h'}+\hat\Sigma_X^{6kh'}}^p \norm{\hat\Sigma_X^{(6k+5)h'}-2\hat\Sigma_X^{(6k+4)h'}+\hat\Sigma_X^{(6k+3)h'}}^p\\
&= \sum_{k=0}^{(6h')^{-1}-1} (h')^{1-p} \bnorm{\Gamma((6kh')\int_{6kh'}^{(6k+3)h'}w_{6k+1}(s)dB'(s)}^p\bnorm{\Gamma(6kh')\int_{(6k+3)h'}^{(6k+6)h'}w_{6k+4}(s)dB'(s)}^p\\
&\quad +{\mathcal O}_{L^1}\Big((h')^{(\beta_b+1/2)\wedge \beta_\Gamma}+(h')^{-1}n^{-1/2}\Big),
\end{align*}
where we only note that the  error bound does not change when pulling $\Gamma(6kh')$ instead of $\Gamma((6k+3)h')$ out of the second stochastic integral.
We  obtain for the conditional expectation by conditional independence of Brownian increments
\begin{align*} &\E\Big[\bnorm{\Gamma(6kh')\int_{6kh'}^{(6k+3)h'}w_{6k+1}(s)dB'(s)}^p\bnorm{\Gamma(6kh')\int_{(6k+3)h'}^{(6k+6)h'}w_{6k+4}(s)dB'(s)}^p\,\Big|\,
{\mathcal F}_{6kh'}\Big]\\
&=\norm{w_k}_{L^2}^{2p}\rho_p(\Gamma(6kh'))^2=(2h')^p\rho_p(\Gamma(6kh'))^2.
\end{align*}
The remaining arguments are then exactly as for the fourth moment result.
\end{proof}

\begin{proof}[Proof of Corollary \ref{Corkappahat2}.]
Set $\underline\lambda_r=\inf_{t\in[0,1]}\lambda_r(\Sigma_X(t))$.
By Proposition \ref{PropCalcCritVal} below (inject $\Sigma_Z=0$, $p=2$, $\underline\lambda_r>0$, $\delta\to 0$ with $\delta^{-1}\le h^{-1/2}$) we have on the event ${\mathcal H}_0^{gap}$
\[ T_{n,h}\le \tfrac{2}{\underline\lambda_r}\Delta_2(\Sigma_X,h)^2+{\mathcal O}_P(\Delta_4(\Sigma_X,h)^2(nh)^{-1/2})=\tfrac{2}{\underline\lambda_r}\Delta_2(\Sigma_X,h)^2+{\mathcal O}_P(h^{1/2}n^{-1/2}),\]
where we use $\Delta_4(\Sigma_X,h)\le h^{1/2}\abs{\Sigma_X}_{B^{1/2}_{4,\infty}}$ and $\Sigma_X\in B^{1/2}_{4,\infty}$ almost surely. Therefore Proposition \ref{PropNQV} and $h+(nh)^{-1}=o(h')$ yield on ${\mathcal H}_0^{gap}$
\[ T_{n,h}\le \tfrac{1}{3\underline\lambda_r}hNV^{(2)}+{\mathcal O}_P(h^{3/2}+h^{1/2}n^{-1/2})= \tfrac{1}{3\underline\lambda_r}h\big(NV^{(2)}+{ o}_P((h')^{1/2})\big).\]
From \eqref{Eqlambdahat} and $h'n^{1/3}\to\infty$ we deduce that then also
\[ T_{n,h}\le  \tfrac{1}{3\hat{\underline\lambda}_r}h\big(NV^{(2)}+{ o}_P((h')^{1/2})\big)\]
holds. By the Central Limit Theorem \ref{ThmNQVCLT}, the consistency results of Proposition \ref{PropVarhat} and Slutsky's lemma we conclude
\[ \Big(3h' \big(\widehat{NV}_{h',n}^{(4)}-\widehat{BNV}^{(2)}_{h',n}\big)\Big)^{-1/2}\Big(\widehat{NV}^{(2)}_{h',n}-NV^{(2)}\Big) \xrightarrow{d} N(0,1).
\]
Another application of Slutsky's lemma thus implies
\[ \liminf_{n\to\infty}\PP\Big(\Big\{\Big(3h' \big(\widehat{NV}_{h',n}^{(4)}-\widehat{BNV}^{(2)}_{h',n}\big)\Big)^{-1/2} \Big(\tfrac{3\hat{\underline\lambda}_r}{h}T_{n,h}-\widehat{NV}^{(2)}_{h',n}\Big)\le q_{1-\alpha;N(0,1)}\Big\}\cap{\mathcal H}_0^{gap}\Big)\]
is at least $1-\alpha$. This gives the result for the test $\phi_{\alpha,n,h,h'}$.
\end{proof}

\end{appendix}

\bibliographystyle{apalike2}
\bibliography{HF-RankTest}

\end{document}